\documentclass[12pt]{amsart}

\title{Total acyclicity of complexes over group algebras}
\author{Ioannis Emmanouil and Olympia Talelli}

\newtheorem{Lemma}{Lemma}[section]
\newtheorem{Proposition}[Lemma]{Proposition}
\newtheorem{Theorem}[Lemma]{Theorem}
\newtheorem{Corollary}[Lemma]{Corollary}

\begin{document}

\begin{abstract}
In this paper, we study group algebras over which modules have a
controlled behaviour with respect to the notions of Gorenstein
homological algebra, namely: (a) Gorenstein projective modules
are Gorenstein flat, (b) any module whose dual is Gorenstein
injective is necessarily Gorentein flat, (c) the Gorenstein
projective cotorsion pair is complete and (d) any acyclic complex
of projective, injective or flat modules is totally acyclic (in
the respective sense). We consider a certain class of groups
satisfying all of these properties and show that it is closed
under the operation ${\scriptstyle{\bf LH}}$ defined by Kropholler
and the operation $\Phi$ defined by the second author. We thus
generalize all previously known results regarding these properties
over group algebras and place these results in an appropriate
framework.
\end{abstract}

\maketitle
\tableofcontents

\addtocounter{section}{-1}
\section{Introduction}

\noindent
The notions of Gorenstein projective, Gorenstein injective and
Gorenstein flat modules generalize the classical notions of
projective, injective and flat modules, respectively. They
were introduced in \cite{EJ1,EJT}, but their origin can be
traced back to earlier work by Auslander and Bridger \cite{A,AB}.
Gorenstein homological algebra, i.e.\ the relative homological
theory based upon these classes of modules, has developed rapidly
during the past decades and found interesting applications in a
wide range of research areas, including representation theory,
the theory of singularities and cohomological group theory.

It is not known whether the relation between Gorenstein projective,
Gorenstein injective and Gorenstein flat modules over a ring $R$ is
aligned to the relation between their homological algebra counterparts.
First of all, it is not known whether

\medskip

\noindent
(a) Gorenstein projective modules are Gorenstein flat.\footnote{As
shown in \cite{H}, this is the case if the ring $R$ is right coherent
and has finite left finitistic dimension.}

\medskip

\noindent
The introduction of projectively coresolved Gorenstein flat modules
\cite{SS} provides a class of modules which are both Gorenstein
projective and Gorenstein flat and reduces the problem of showing
that assertion (a) above is true to the problem of showing that any
Gorenstein projective module is projectively coresolved Gorenstein
flat.

\medskip

\noindent
On the other hand, the Pontryagin duals of Gorenstein flat modules
are always Gorenstein injective (right) modules, but it is not known
whether

\medskip

\noindent
(b) any module whose dual is Gorenstein injective is necessarily
Gorenstein flat.\footnote{As shown in \cite{H}, this is the case
if the ring $R$ is right coherent.}

\medskip

\noindent
It is shown in \cite{SS} that the class of Gorenstein flat (resp.\
projectively coresolved Gorenstein flat) modules forms the left
hand side of a complete hereditary cotorsion pair, whereas the
class of Gorenstein injective modules forms the right hand side
of a complete hereditary cotorsion pair. Even though the class
of Gorenstein projective modules forms the left hand side of a
hereditary cotorsion pair as well (cf.\ \cite{CIS}), it is an
open question whether

\medskip

\noindent
(c) the Gorenstein projective cotorsion pair is complete.

\medskip

\noindent
In view of the very definition of the Gorenstein projective (resp.\
Gorenstein flat or Gorenstein injective) modules, these are obtained
by considering certain acyclic complexes of projective (resp.\ flat
or injective) modules, the so-called totally acyclic complexes of
projective (resp.\ flat or injective) modules. A basic problem
regarding these Gorenstein modules is to analyse the relation between
acyclicity and total acyclicity for such complexes. The latter
problem has been studied over commutative Noetherian rings that
admit a dualizing complex in \cite{IK}. In general, one could look
for conditions under which

\medskip

\noindent
(d) any acyclic complex of projective (resp.\ flat or injective)
modules is totally acyclic.

\medskip

\noindent
Some of the four general problems displayed above have been studied
in the literature, in the case where $R=kG$ is the group algebra of
a group $G$ with coefficients in a commutative ring $k$. The special
case where $k$ is the ring of integers and $G$ is a group contained
in Kropholler's class ${\scriptstyle{\bf LH}}\mathfrak{F}$ of
hierarchically decomposable groups with base the class of finite groups
(cf.\ \cite{K}) has been examined in \cite{DT}; it was conjectured
therein that any acyclic complex of projective $\mathbb{Z}G$-modules
is totally acyclic. The results obtained in \cite{DT} have been extended
in \cite{Bi2,MS} to the case where $k$ is a Noetherian ring of finite
global dimension and $G$ is either an
${\scriptstyle{\bf LH}}\mathfrak{F}$-group or a group of type $\Phi$,
as defined by the second author in \cite{T}. The case where the ring
$k$ is Gorenstein regular (resp.\ weakly Gorenstein regular) and $G$
has finite cohomological (resp.\ homological) dimension over $k$ was
examined in \cite{R,RY}; see also \cite{ET2}.

In this paper, we generalize and present a common framework that covers
the above results. Roughly speaking, we show that the classes of groups
for which the four general questions considered above admit a positive
answer are closed under Kropholler's ${\scriptstyle{\bf LH}}$ operation
and the second author's $\Phi$ operation; see Theorems 3.3, 4.2 and 5.2
for precise statements. In order to give a flavour of our results, we
consider the case where $k=\mathbb{Z}$ is the ring of integers. Let
$\mathfrak{Z}$ be the class of those groups $G$, for which the following
two conditions hold:

(1) any acyclic complex of injective $\mathbb{Z}G$-modules is totally
acyclic and

(2) any $\mathbb{Z}G$-module whose dual is Gorenstein injective is
necessarily Gorenstein flat.
\newline
We recall that a group $G$ is of type $\Phi$
if the $\mathbb{Z}G$-modules of finite projective dimension are precisely
those $\mathbb{Z}G$-modules that have finite projective dimension over
$\mathbb{Z}H$ for any finite subgroup $H \subseteq G$. These groups were
introduced in \cite{T}, in view of their importance for the construction
of finite dimensional models for the classifying space for proper actions.
More generally, for any group class $\mathfrak{C}$ we may define
$\Phi_{proj}\mathfrak{C}$ to be the class consisting of those groups $G$,
over which the $\mathbb{Z}G$-modules of finite projective dimension
are precisely those $\mathbb{Z}G$-modules that have uniformly finite
projective dimension over $\mathbb{Z}H$ for any $\mathfrak{C}$-subgroup
$H \subseteq G$. As a variant of that definition, for any group class
$\mathfrak{C}$ we may consider the class $\Phi_{inj} \mathfrak{C}$
consisting of those groups $G$, over which the $\mathbb{Z}G$-modules of
finite injective dimension are precisely those $\mathbb{Z}G$-modules that
have uniformly finite injective dimension over $\mathbb{Z}H$ for any
$\mathfrak{C}$-subgroup $H \subseteq G$. The following result is a
consequence of Proposition 5.1, Theorem 5.2 and Proposition 5.5.

\vspace{0.1in}

\noindent
{\bf Theorem.}
{\em Let $\mathfrak{Z}$ be the class of groups defined above. Then:

(i) The class $\mathfrak{Z}$ contains all groups $G$ for which
$\mathbb{Z}G$ is weakly Gorenstein regular. In particular,
$\mathfrak{Z}$ contains all groups of type $\Phi$ and hence all
finite groups.

(ii) ${\scriptstyle{\bf LH}}\mathfrak{Z} = \mathfrak{Z}$ and
$\Phi_{inj}\mathfrak{Z} = \mathfrak{Z}$.

(iii) If $G$ is a $\mathfrak{Z}$-group, then the following hold:
(3) any Gorenstein projective $\mathbb{Z}G$-module is Gorenstein
flat, (4) the Gorenstein projective cotorsion pair is complete and
(5) any acyclic complex of projective (resp.\ flat)
$\mathbb{Z}G$-modules is totally acyclic.}

\vspace{0.1in}

\noindent
The point here is that property (1) in the definition of the class
$\mathfrak{Z}$ implies properties (3), (4) and (5) in assertion
(iii) of the theorem. This is proved using an interesting relation,
which is universally valid over any ring and not only over group
algebras, between modules that appear as kernels of acyclic complexes
of projective, flat or injective modules; this relation is described
in Section 2 of the paper. On the other hand, if $\mathfrak{Z}_0$ is
any subclass of $\mathfrak{Z}$, then (ii) implies that $\mathfrak{Z}$
also contains the smallest class of groups, say
$\overline{\mathfrak{Z}_0}$, that contains $\mathfrak{Z}_0$ and is
both ${\scriptstyle{\bf LH}}$-closed and $\Phi_{inj}$-closed. As we
explain in $\S $1.IV, the groups in $\overline{\mathfrak{Z}_0}$ admit
a hierarchical description with base the group class $\mathfrak{Z}_0$,
\`{a} la Kropholler. In principle, this provides us with an abundance
of $\mathfrak{Z}$-groups. For example, groups in the
${\scriptstyle{\bf LH}}$-closure of the particular subclass of
$\mathfrak{Z}$ provided in (i) include all cases that were
previously considered in the literature, regarding only some of the
four problems considered here. The reader may consult the statements
of the results in Sections 3, 4 and 5 for more general choices of the
coefficient ring $k$.

As an addendum to the above result, and more precisely to property
(4) in assertion (iii) therein, we note that an explicit set of modules
that cogenerate the Gorenstein projective cotorsion pair is obtained in
\cite{Ke}, if $G$ is in Kropholler's class of hierarchically decomposable
groups with base the class of those groups for which the associated group
algebra is Gorenstein regular. Any cotorsion pair cogenerated by a set is
known to be complete; cf.\ \cite{ET}. We also note that property (5) in
assertion (iii) implies that the tensor product $M \otimes_{\mathbb{Z}} N$
of two $\mathbb{Z}G$-modules, endowed with the diagonal action of $G$, is
Gorenstein projective (resp.\ Gorenstein flat) if $M$ is Gorenstein
projective (resp.\ Gorenstein flat) and $N$ is $\mathbb{Z}$-free (resp.\
$\mathbb{Z}$-flat). On the other hand, property (1) in the definition of
$\mathfrak{Z}$ implies that the diagonal $\mathbb{Z}G$-module
$\mbox{Hom}_{\mathbb{Z}}(M,N)$ is Gorenstein injective, provided that $N$
is Gorenstein injective and $M$ is $\mathbb{Z}$-flat. Moreover, for any
subgroup $H \subseteq G$, the validity of property (5) for $H$ implies
that the restriction of any Gorenstein projective (resp.\ Gorenstein
flat) $\mathbb{Z}G$-module is a Gorenstein projective (resp.\ Gorenstein
flat) $\mathbb{Z}H$-module. Analogously, if $H$ has property (1) in the
definition of $\mathfrak{Z}$, then the restriction of any Gorenstein
injective $\mathbb{Z}G$-module is a Gorenstein injective $\mathbb{Z}H$-module.

The contents of the paper are as follows: The preliminary Section
1 collects the background notions and results that are needed in
the paper. In Section 2, we study the relation between modules
(over a general ring) that appear as cokernels of acyclic complexes
of flat modules to those modules that appear as cokernels of acyclic
complexes of projective modules, as well as to those modules that
appear as kernels of acyclic complexes of injective modules. Then,
in Section 3, we examine groups $G$ for which the cokernels of any
acyclic complex of projective (resp.\ flat) modules over the group
algebra $kG$ are Gorenstein projective (resp.\ Gorenstein flat) and
show that the class of these groups is closed under the operations
${\scriptstyle{\bf LH}}$ and $\Phi$. The analogous property regarding
the kernels of acyclic complexes of injective $kG$-modules is examined
in Section 4. In Section 5, we consider a class of groups that
have all properties examined in the discussion above (in the form of
assertions (a), (b), (c) and (d)) and show that this class is closed
under the operations ${\scriptstyle{\bf LH}}$ and $\Phi_{inj}$.

\vspace{0.1in}

\noindent
{\em Notations and terminology.}
All rings are unital and associative. If $R$ is a ring, we let
${\tt Proj}(R)$, ${\tt Inj}(R)$, ${\tt Flat}(R)$ and ${\tt Cotor}(R)$
be the classes of projective, injective, flat and cotorsion left
$R$-modules, respectively. If $M$ is a left $R$-module, its
Pontryagin dual is the right $R$-module
$DM = \mbox{Hom}_{\mathbb{Z}}(M,\mathbb{Q}/\mathbb{Z})$ consisting of
all additive maps from $M$ to the abelian group $\mathbb{Q}/\mathbb{Z}$.

\section{Preliminaries}

\noindent
In this Section, we record the prerequisite notions and results
that will be used in the sequel. These notions concern the basics
of Gorenstein homological algebra, special features regarding
modules over group algebras and the closure of group classes
under certain operations.

\vspace{0.1in}

\noindent
{\sc I.\ Gorenstein modules.}
We fix a ring $R$ and consider, unless otherwise specified, only
left $R$-modules. An acyclic complex of projective modules is
called totally acyclic if it remains acyclic after applying the
functor $\mbox{Hom}_R(\_\!\_,P)$ for any projective module $P$.
A module $M$ is Gorenstein projective if it is a cokernel of a
totally acyclic chain complex of projective modules; we denote by
${\tt GProj}(R)$ the class of these modules. Dually, an acyclic
complex of injective modules is called totally acyclic if it
remains acyclic after applying the functor $\mbox{Hom}_R(I,\_\!\_)$
for any injective module $I$. A module $N$ is Gorenstein injective
if it is a kernel of a totally acyclic cochain complex of injective
modules; we denote by ${\tt GInj}(R)$ the class of these modules.
An acyclic complex of flat modules is called totally acyclic if
it remains acyclic after applying the functor $J \otimes_R\_\!\_$
for any injective right module $J$. A module $L$ is Gorenstein
flat if it is a cokernel of a totally acyclic chain complex of
flat modules; we denote by ${\tt GFlat}(R)$ the class of these
modules. The standard reference for these notions is Holm's paper
\cite{H}. A special class of Gorenstein flat modules was introduced
in \cite{SS}: We say that a module $K$ is projectively coresolved
Gorenstein flat if it is a cokernel of an acyclic chain complex
of projective modules which remains acyclic after applying the
functor $J \otimes_R\_\!\_$ for any injective right module $J$;
we denote by ${\tt PGF}(R)$ the class of these modules. As shown
in \cite[Theorem 4.4]{SS}, any projectively coresolved Gorenstein
flat module is Gorenstein projective, i.e.\
${\tt PGF}(R) \subseteq {\tt GProj}(R)$. We also note that the
Pontryagin dual $DL$ of any Gorenstein flat module $L$ is Gorenstein
injective (as a right module), i.e.\
$D{\tt GFlat}(R) \subseteq {\tt GInj}(R^{op})$; this is proved in
\cite[Theorem 3.6]{H}.

The Gorenstein projective dimension $\mbox{Gpd}_RM$ of a module $M$
is the length of a shortest resolution of $M$ by Gorenstein projective
modules. If there is no such resolution of finite length, then we write
$\mbox{Gpd}_RM = \infty$. Analogously, we may define the (projectively
coresolved) Gorenstein flat dimension of modules. Dually, the Gorenstein
injective dimension $\mbox{Gid}_RN$ of a module $N$ is the length of a
shortest coresolution of $N$ by Gorenstein injective modules. If there
is no such resolution of finite length, then we write
$\mbox{Gid}_RN = \infty$. Regarding the duality relation between
Gorenstein flat and Gorenstein injective modules, we record the
following criterion for a module to be Gorenstein flat; it is due
to Bouchiba \cite[Theorem 4(5)]{Bou}.

\begin{Proposition}
The following conditions are equivalent for a module $M$:

(i) $M \in {\tt GFlat}(R)$.

(ii) $M$ has finite Gorenstein flat dimension and
     $DM \in {\tt GInj}(R^{op})$. \hfill $\Box$
\end{Proposition}

\noindent
As shown by Beligiannis \cite[$\S $6]{Bel}, the finiteness of the
Gorenstein projective dimension $\mbox{Gpd}_RM$ of a module $M$ is
equivalent to the existence of a complete projective resolution of
$M$, i.e.\ to the existence of a totally acyclic complex of projective
modules that coincides in sufficiently large degrees with an ordinary
projective resolution of $M$. The existence of complete projective
resolutions of modules has been studied by Gedrich and Gruenberg
\cite{GG}, as well as by Cornick and Kropholler \cite{CK}. A
characterization of the rings $R$ over which all modules admit
complete projective resolutions is obtained in \cite{CK}, in terms
of the finiteness of the two invariants $\mbox{spli} \, R$ and
$\mbox{silp} \, R$. Here, $\mbox{spli} \, R$ is the supremum of the
projective dimensions of injective modules and $\mbox{silp} \, R$
is the supremum of the injective dimensions of projective modules.
That result by Cornick and Kropholler was alternatively proved using
the language of Gorenstein projective dimension by Bennis and Mahbou
in \cite{BM}, where the notion of the Gorenstein global dimension of
the ring was introduced. More precisely, the (left) Gorenstein global
dimension $\mbox{Ggl.dim} \, R$ of the ring $R$ is defined as the
supremum of the Gorenstein projective dimensions of all modules or,
equivalently, as the supremum of the Gorenstein injective dimensions
of all modules. Its finiteness is equivalent to the finiteness of both
$\mbox{spli} \, R$ and $\mbox{silp} \, R$, in which case we have
\[ \mbox{Ggl.dim} \, R =  \mbox{spli} \, R = \mbox{silp} \, R <
   \infty . \]
Following Beligiannis \cite{Bel}, we say that $R$ is (left)
Gorenstein regular if $\mbox{Ggl.dim} \, R < \infty$.

The characterization of the finiteness of the Gorenstein weak
global dimension $\mbox{Gwgl.dim} \, R$ of $R$, which is defined
as the supremum of the Gorenstein flat dimensions of all modules,
turned out to be more difficult to obtain. Here, the homological
invariants that play a role are $\mbox{sfli} \, R$, the supremum
of the flat dimensions of injective modules, and its analogue
$\mbox{sfli} \, R^{op}$ for the opposite ring $R^{op}$. It has
been shown by Christensen, Estrada and Thompson in \cite{CET}
that $\mbox{Gwgl.dim} \, R < \infty$ if and only if both
$\mbox{sfli} \, R$ and $\mbox{sfli} \, R^{op}$ are finite, in
which case we have
\[ \mbox{Gwgl.dim} \, R =  \mbox{sfli} \, R = \mbox{sfli} \, R^{op}
   < \infty . \]
We note that the above result admits a simpler proof if $R$ is
a commutative ring or, more generally, if $R$ is isomorphic with
its opposite $R^{op}$; cf.\ \cite{E1}. We say that $R$ is weakly
Gorenstein regular if $\mbox{Gwgl.dim} \, R < \infty$.

\vspace{0.1in}

\noindent
{\sc II.\ Cotorsion pairs.}
We work over a ring $R$ and consider the $\mbox{Ext}^1$-pairing,
which induces an orthogonality relation between modules. If {\tt S}
is a class of modules, then the left orthogonal $^{\perp}{\tt S}$
of ${\tt S}$ is the class consisting of those modules $M$, which
are such that $\mbox{Ext}^1_R(M,S) = 0$ for all $S \in {\tt S}$.
Analogously, the right orthogonal ${\tt S}^{\perp}$ of {\tt S} is
the class consisting of those modules $N$, which are such that
$\mbox{Ext}^1_R(S,N) = 0$ for all $S \in {\tt S}$. If
${\tt U},{\tt V}$ are two module classes, we say that the pair
$({\tt U},{\tt V})$ is a cotorsion pair (cf.\ \cite{EJ2}) if
${\tt U} = {^{\perp}}{\tt V}$ and ${\tt U}^{\perp} = {\tt V}$.
The cotorsion pair is hereditary if $\mbox{Ext}^i_R(U,V)=0$ for
all $i>0$ and all $U \in {\tt U}$ and $V \in {\tt V}$; it is
complete if for any module $M$ there exist short exact sequences
\[ 0 \longrightarrow V \longrightarrow U \longrightarrow M
     \longrightarrow 0
   \;\;\; \mbox{and } \;\;\;
   0 \longrightarrow M \longrightarrow V' \longrightarrow U'
     \longrightarrow 0 , \]
where $U,U' \in {\tt U}$ and $V,V' \in {\tt V}$.

Gorenstein modules provide non-trivial examples of cotorsion pairs:
It is proved by \v{S}aroch and \v{S}t$\!$'ov\'{i}\v{c}ek \cite{SS}
that $\left( {\tt PGF}(R),{\tt PGF}(R)^{\perp} \right)$,
$\left( {\tt GFlat}(R),{\tt GFlat}(R)^{\perp} \right)$ and
$\left( ^{\perp}{\tt GInj}(R),{\tt GInj}(R) \right)$
are all complete and hereditary cotorsion pairs. We note that
${\tt PGF}(R) \cap {\tt PGF}(R)^{\perp} = {\tt Proj}(R)$,
${\tt GFlat}(R) \cap {\tt GFlat}(R)^{\perp} =
 {\tt Flat}(R) \cap {\tt Cotor}(R)$
and $^{\perp}{\tt GInj}(R) \cap {\tt GInj}(R) = {\tt Inj}(R)$.
The situation for Gorenstein projective modules seems to be
more subtle: It is shown in \cite{CIS} that
$\left( {\tt GProj}(R),{\tt GProj}(R)^{\perp} \right)$ is a
hereditary cotorsion pair as well, but its completeness appears
to depend upon set-theoretic assumptions.

\vspace{0.1in}

\noindent
{\sc III.\ Modules over group rings.}
We consider a commutative ring $k$ and let $G$ be a group. Using
the diagonal action of $G$, the tensor product $M \otimes_k N$ of
two $kG$-modules $M,N$ is also a $kG$-module; we define
$g \cdot (x \otimes y) = gx \otimes gy \in M \otimes_k N$ for any
$g \in G$, $x \in M$ and $y \in N$. If $M$ is a flat $kG$-module
and $N$ is flat as a $k$-module, then the diagonal $kG$-module
$M \otimes_k N$ is flat as well. Indeed, since any flat $kG$-module
is the filtered colimit of free $kG$-modules and the class of flat
$kG$-modules is closed under filtered colimits and direct sums, we
may reduce the proof of this claim to the special case where $M=kG$.
In that case, the result follows from
\cite[Chapter III, Corollary 5.7]{Br}, which implies that the
diagonal $kG$-module $kG \otimes_kN$ is isomorphic with the induced
$kG$-module $\mbox{ind}_1^G\mbox{res}_1^GN$.

In Section 5, we shall use assertion (iii) of the following result.
Even though that assertion is perhaps folklore knowledge (see
\cite[Example 3.8]{CR} for the case of Gorenstein projective
modules), we provide a detailed proof using the technique
introduced by Cornick and Kropholler in \cite[Theorem 3.5]{CK}.

\begin{Proposition}
Let $k$ be a commutative ring and $G$ a finite group.

(i) An acyclic complex of flat $kG$-modules $F$ is totally
acyclic if and only if the acyclic complex of flat $k$-modules
$\mbox{res}_1^GF$ is totally acyclic.

(ii) If $M,N$ are $kG$-modules with $M$ Gorenstein flat and $N$
$k$-flat, then the diagonal $kG$-module $M \otimes_kN$ is also
Gorenstein flat.

(iii) A $kG$-module $M$ is Gorenstein flat if and only if the
induced $kG$-module $\mbox{ind}_1^G\mbox{res}_1^GM$ is Gorenstein
flat.
\end{Proposition}
\vspace{-0.05in}
\noindent
{\em Proof.}
(i) This is clear, since any injective $kG$-module $I$ is a
direct summand of a $kG$-module of the form
$\mbox{coind}_1^GJ \simeq \mbox{ind}_1^GJ$, where
$J \in {\tt Inj}(k)$, and
$F \otimes_{kG} \mbox{ind}_1^GJ \simeq \mbox{res}_1^GF \otimes_k J$.

(ii) There exists a totally acyclic complex of flat $kG$-modules
$F$, such that $M = C_0F$. Since $N$ is $k$-flat, the complex of
diagonal $kG$-modules $F \otimes_kN$ is acyclic, consists of flat
$kG$-modules and $M \otimes_kN = C_0(F \otimes_kN)$. In order to
show that $F \otimes_kN$ is totally acyclic (as a complex of flat
$kG$-modules), we invoke (i) and reduce to showing that the complex
$F \otimes_kN \otimes_kJ$ is acyclic for any injective $k$-module
$J$. The latter claim follows since the complex $F \otimes_kJ$ is
acyclic (using again (i) above) and $N$ is $k$-flat.

(iii) As we have noted above, \cite[Chapter III, Corollary 5.7]{Br}
implies that the induced $kG$-module $\mbox{ind}_1^G\mbox{res}_1^GM$
is isomorphic with the diagonal $kG$-module $M \otimes_kkG$. If $M$
is Gorenstein flat, then $M \otimes_kkG$ is Gorenstein flat as well,
in view of (ii) above. Conversely, assume that the diagonal
$kG$-module $M \otimes_kkG$ is Gorenstein flat and consider the
$k$-split short exact sequences of $kG$-modules
\[ 0 \longrightarrow L \longrightarrow kG
   \stackrel{\varepsilon}{\longrightarrow} k \longrightarrow 0
   \;\;\; \mbox{and} \;\;\;
   0 \longrightarrow k \stackrel{\eta}{\longrightarrow} kG
   \longrightarrow N \longrightarrow 0 , \]
where $\varepsilon$ is the augmentation map and $\eta$ maps
$1$ onto $\sum_{g \in G}g$. We note that the $kG$-modules $L,N$
are both $k$-free. Let $L^{\otimes i}$ and $N^{\otimes i}$ denote
the $i$-th tensor powers over $k$ of $L$ and $N$ respectively for
all $i \geq 0$. Tensoring the former of these short exact sequences
with $M \otimes_k L^{\otimes i}$, $i \geq 0$, and the latter sequence
with $M \otimes_k N^{\otimes i}$, $i \geq 0$, we obtain two sequences
of short exact sequences that may be spliced together to a doubly
infinite acyclic complex $F$ of $kG$-modules
\[ \cdots \! \longrightarrow F_2 \longrightarrow F_1
   \longrightarrow F_0 \longrightarrow F^0 \longrightarrow F^1
   \longrightarrow F^2 \longrightarrow \cdots . \]
Here, $F_i = M \otimes_kkG \otimes_k L^{\otimes i}$ and
$F^i = M \otimes_kkG \otimes_k N^{\otimes i}$ for all $i \geq 0$.
In view of our assumption and (ii) above, it follows that $F$
is a complex of Gorenstein flat $kG$-modules. The cokernels of
$F$ are the $kG$-modules $(M \otimes_k L^{\otimes i})_{i>0}$,
$M$ and $(M \otimes_k N^{\otimes i})_{i>0}$.

We shall prove that the complex $I \otimes_{kG} F$ is acyclic for
any injective $kG$-module $I$; invoking \cite{BK}, it will then
follow that the $kG$-module $M$ is indeed Gorenstein flat. To that
end, we fix an injective $kG$-module $I$ and note that the $k$-split
monomorphism $\eta$ induces a monomorphism of $kG$-modules with
diagonal action
\[ 1 \otimes \eta : I = I \otimes_kk \longrightarrow I \otimes_kkG . \]
This monomorphism is necessarily split, since $I$ is an injective
$kG$-module. As the $kG$-module $I$ is a direct summand of $I \otimes_kkG$,
the acyclicity of the complex $I \otimes_{kG} F$ will follow if we
show the acyclicity of the complex
$(I \otimes_kkG) \otimes_{kG} F \simeq I \otimes_{kG} (kG \otimes_kF)$.\footnote{The
isomorphism follows since both sides are quotients of the tensor product
$(I \otimes_kkG) \otimes_kF \simeq I \otimes_k (kG \otimes_k F)$ modulo
the relations of the form $gx \otimes gy \otimes gz = x \otimes y \otimes z$,
where $g \in G$, $x \in I$, $y \in kG$, $z \in F$.} The complex of
$kG$-modules $kG \otimes_k F$ (with diagonal action) is acyclic and
has cokernels $(M \otimes_kkG \otimes_k L^{\otimes i})_{i>0}$,
$M \otimes_kkG$ and $(M  \otimes_kkG \otimes_k N^{\otimes i})_{i>0}$.
Since all of these $kG$-modules are Gorenstein flat, the complex
$kG \otimes_kF$ decomposes into short exact sequences with all
three terms Gorenstein flat. The functor $\mbox{Tor}_1^{kG}(I,\_\!\_)$
vanishes on Gorenstein flat modules and hence the latter short exact
sequences remain exact after applying the functor $I \otimes_{kG} \_\!\_\,$.
The complex $I \otimes_{kG} (kG \otimes_kF)$ is thus acyclic, as needed.
\hfill $\Box$

\vspace{0.1in}

\noindent
Let $H \subseteq G$ be a subgroup. The reader is referred to \cite{Br}
for the basic properties of the restriction functor $\mbox{res}_H^G$
from $kG$-modules to $kH$-modules, as well as the induction and
coinduction functors $\mbox{ind}_H^G$ and $\mbox{coind}_H^G$ from
$kH$-modules to $kG$-modules. In the sequel, we shall make use of
the following four properties regarding these functors. Properties
(a) and (b) are standard, whereas properties (c) and (d) are simple
consequences of the definitions.
\newline
(a) Any $kH$-module $M$ is contained as a direct summand in
$\mbox{res}_H^G\mbox{ind}_H^GM$ and $\mbox{res}_H^G\mbox{coind}_H^GM$.
\newline
(b) If $H$ has finite index in $G$, then the functors $\mbox{ind}_H^G$
and $\mbox{coind}_H^G$ are naturally isomorphic.
\newline
(c) There is a natural isomorphism of functors
$D \circ \mbox{ind}_H^G \simeq \mbox{coind}_H^G \circ D$, where $D$
denotes the Pontryagin duality functor in the respective module
categories.
\newline
(d) The induction functor $\mbox{ind}_H^G$ maps ${\tt PGF}(kH)$
into ${\tt PGF}(kG)$ and ${\tt GFlat}(kH)$ into ${\tt GFlat}(kG)$,
whereas the coinduction functor $\mbox{coind}_H^G$ maps
${\tt GInj}(kH)$ into ${\tt GInj}(kG)$.

\vspace{0.1in}

\noindent
{\sc IV.\ Group class operations.}
An operation $\mathcal{O}$ on group classes assigns a group class
$\mathcal{O}\mathfrak{C}$ to each group class $\mathfrak{C}$, in
such a way that the following two properties are satisfied:
\newline
(i) $\mathfrak{C} \subseteq \mathcal{O}\mathfrak{C}$ for any group
class $\mathfrak{C}$ and
\newline
(ii) $\mathcal{O}\mathfrak{C}_1 \subseteq \mathcal{O}\mathfrak{C}_2$
for any two group classes $\mathfrak{C}_1,\mathfrak{C}_2$ with
$\mathfrak{C}_1 \subseteq \mathfrak{C}_2$.
\newline
There is a notion of continuity for operations on group classes. In
order to define that notion, assume that $\mathfrak{C}_{\alpha}$ is
a group class defined for each ordinal number $\alpha$ and
$\mathfrak{C}_{\alpha} \subseteq \mathfrak{C}_{\beta}$ for all
$\alpha \leq \beta$. Then, we may consider the group class
$\mathfrak{C}$ consisting of those groups which are contained
in $\mathfrak{C}_{\alpha}$ for some ordinal number $\alpha$. We
say that an operation $\mathcal{O}$ is continuous if it has the
following additional property:
\newline
(iii) Whenever $\mathfrak{C}_{\alpha}$, $\alpha$ an ordinal, and
$\mathfrak{C}$ are group classes as above, the group class
$\mathcal{O}\mathfrak{C}$ is the class consisting of those groups
which are contained in $\mathcal{O}\mathfrak{C}_{\alpha}$ for some
$\alpha$.

We say that a group class $\mathfrak{C}$ is closed under an
operation $\mathcal{O}$ if $\mathcal{O}\mathfrak{C} = \mathfrak{C}$.
The $\mathcal{O}$-closure of a group class $\mathfrak{C}$ is
the smallest $\mathcal{O}$-closed class of groups containing
$\mathfrak{C}$; it is precisely the class consisting of those
groups $G$ which are contained in $\mathfrak{D}$ for any
$\mathcal{O}$-closed class $\mathfrak{D}$ containing $\mathfrak{C}$.
There is a hierachical description of the $\mathcal{O}$-closure of a
group class $\mathfrak{C}$, in the special case where $\mathcal{O}$
is a continuous operation. More generally, let
$\mathcal{O}_1, \mathcal{O}_2, \ldots , \mathcal{O}_n$ be continuous
operations on group classes and consider a group class $\mathfrak{C}$.
Then, we may define the classes $\mathfrak{C}_{\alpha}$ for any ordinal
number $\alpha$, using transfinite induction, as follows:
\newline
(a) $\mathfrak{C}_0 = \mathfrak{C}$,
\newline
(b) $\mathfrak{C}_{\alpha +1} =
     \mathcal{O}_1\mathfrak{C}_{\alpha} \cup
     \mathcal{O}_2\mathfrak{C}_{\alpha} \cup \ldots \cup
     \mathcal{O}_n\mathfrak{C}_{\alpha}$
     for any ordinal $\alpha$ and
\newline
(c) $\mathfrak{C}_{\alpha} =
     \bigcup_{\beta < \alpha}\mathfrak{C}_{\beta}$
    for any limit ordinal $\alpha$.
\newline
Let $\overline{\mathfrak{C}}$ be the class consisting of those
groups which are contained in $\mathfrak{C}_{\alpha}$, for some
ordinal $\alpha$. It is easily seen that $\overline{\mathfrak{C}}$
is the smallest class of groups that contains $\mathfrak{C}$ and
is closed under all of the operations
$\mathcal{O}_1, \mathcal{O}_2, \ldots , \mathcal{O}_n$.

We are interested in the operation on group classes that was
introduced by Kropholler in \cite{K}. If $\mathfrak{C}$ is
any group class, then we may define for each ordinal $\alpha$
a group class ${\scriptstyle{{\bf H}}}_{\alpha}\mathfrak{C}$,
using transfinite induction, as follows: First of all, we
let ${\scriptstyle{{\bf H}}}_0\mathfrak{C} = \mathfrak{C}$.
For any ordinal $\alpha >0$, the class
${\scriptstyle{{\bf H}}}_{\alpha}\mathfrak{C}$ consists of those
groups $G$ which admit a cellular action on a finite dimensional
contractible CW-complex $X$, in such a way that the isotropy subgroup
of each cell belongs to ${\scriptstyle{{\bf H}}}_{\beta}\mathfrak{C}$
for some ordinal $\beta < \alpha$ (that may depend upon the cell).
If $G,X$ are as above, $k$ is a commutative ring and $M$ is a
$kG$-module, then the cellular chain complex of $X$ induces two
exact sequences of $kG$-modules
\[ 0 \longrightarrow M_r \longrightarrow \cdots
     \longrightarrow M_1 \longrightarrow M_0
     \longrightarrow M \longrightarrow 0 \]
and
\[ 0 \longrightarrow M \longrightarrow M^0
     \longrightarrow M^1 \longrightarrow \cdots
     \longrightarrow M^r \longrightarrow 0 , \]
which are often useful. Here, the $M_i$'s (resp.\ the $M^i$'s) are
direct sums (resp.\ direct products) of $kG$-modules of the form
$\mbox{ind}_H^G\mbox{res}_H^GM$ (resp.\
$\mbox{coind}_H^G\mbox{res}_H^GM$), where $H$ is an
${\scriptstyle{{\bf H}}}_{\beta}\mathfrak{C}$-subgroup of $G$ for
some $\beta < \alpha$. We say that a group belongs to
${\scriptstyle{{\bf H}}}\mathfrak{C}$ if it belongs to
${\scriptstyle{{\bf H}}}_{\alpha}\mathfrak{C}$ for some ordinal
$\alpha$. Then, the class ${\scriptstyle{{\bf LH}}}\mathfrak{C}$
consists of those groups $G$, all of whose finitely generated
subgroups $H$ are contained in some
${\scriptstyle{{\bf H}}}\mathfrak{C}$-subgroup $K=K(H) \subseteq G$.

In this paper, we are primarily interested in the special case
of subgroup-closed group classes. If $\mathfrak{C}$ is such a
subgroup-closed class, then we may use transfinite induction and
show that the class ${\scriptstyle{{\bf H}}}_{\alpha}\mathfrak{C}$
is subgroup-closed for any ordinal $\alpha$. In particular, the
class ${\scriptstyle{{\bf H}}}\mathfrak{C}$ is subgroup-closed
as well. It follows that ${\scriptstyle{{\bf LH}}}\mathfrak{C}$
is also subgroup-closed; in fact, in this case,
${\scriptstyle{{\bf LH}}}\mathfrak{C}$ consists of those groups
all of whose finitely generated subgroups are
${\scriptstyle{{\bf H}}}\mathfrak{C}$-groups.

\begin{Lemma}
(i) The operation ${\scriptstyle{{\bf H}}}_{\alpha}$ is continuous
for any ordinal $\alpha$.

(ii) The operation ${\scriptstyle{{\bf H}}}$ is continuous.

(iii) The operation ${\scriptstyle{{\bf LH}}}$ is continuous.
\end{Lemma}
\vspace{-0.05in}
\noindent
{\em Proof.}
Assume that $\mathfrak{C}_{\beta}$ is a group class defined for
any ordinal $\beta$, such that
$\mathfrak{C}_{\beta} \subseteq \mathfrak{C}_{\beta'}$ whenever
$\beta \leq \beta'$. We also consider the group class $\mathfrak{C}$
consisting of those groups which are contained in
$\mathfrak{C}_{\beta}$ for some ordinal $\beta$.

(i) We use transfinite induction on $\alpha$ to show that any
group $G \in {\scriptstyle{{\bf H}}}_{\alpha}\mathfrak{C}$ is
contained in
${\scriptstyle{{\bf H}}}_{\alpha}\mathfrak{C}_{\beta}$ for a
suitable ordinal $\beta$. This is clear if $\alpha = 0$, since
${\scriptstyle{{\bf H}}}_0\mathfrak{C} = \mathfrak{C}$ and
${\scriptstyle{{\bf H}}}_0\mathfrak{C}_{\beta} =
 \mathfrak{C}_{\beta}$
for all ordinals $\beta$. We now assume that $\alpha >0$ is an
ordinal and the result holds for all ordinals $\alpha' < \alpha$.
If $G \in {\scriptstyle{{\bf H}}}_{\alpha}\mathfrak{C}$, then
$G$ admits a cellular action on a finite dimensional contractible
CW-complex, in such a way that the isotropy subgroup of each cell
$c$ belongs to ${\scriptstyle{{\bf H}}}_{\alpha'(c)}\mathfrak{C}$
for some ordinal $\alpha'(c) < \alpha$. Using the inductive assumption
and the fact that the subgroups of $G$ form a {\em set}, it follows
that there exists an ordinal $\beta$, such that the isotropy subgroup
of any cell $c$ belongs to
${\scriptstyle{{\bf H}}}_{\alpha'(c)}\mathfrak{C}_{\beta}$ for some
ordinal $\alpha'(c) < \alpha$. Then,
$G \in {\scriptstyle{{\bf H}}}_{\alpha}\mathfrak{C}_{\beta}$.

(ii) Let $G$ be an ${\scriptstyle{{\bf H}}}\mathfrak{C}$-group.
Then, $G \in {\scriptstyle{{\bf H}}}_{\alpha}\mathfrak{C}$ for
some ordinal $\alpha$. Since the operation
${\scriptstyle{{\bf H}}}_{\alpha}$ is continuous (in view of
(i) above), we conclude that
$G \in {\scriptstyle{{\bf H}}}_{\alpha}\mathfrak{C}_{\beta}$
for some ordinal $\beta$. It follows that
$G \in {\scriptstyle{{\bf H}}}\mathfrak{C}_{\beta}$.

(iii) Let $G$ be an ${\scriptstyle{{\bf LH}}}\mathfrak{C}$-group.
If $H \subseteq G$ is a finitely generated subgroup, then there
exists an ${\scriptstyle{{\bf H}}}\mathfrak{C}$-subgroup
$K(H) \subseteq G$ containing $H$. Since the operation
${\scriptstyle{{\bf H}}}$ is continuous (in view of (ii)
above) and the finitely generated subgroups of $G$ form a
{\em set}, we can find a suitable ordinal $\beta$, such
that
$K(H) \in {\scriptstyle{{\bf H}}}\mathfrak{C}_{\beta}$
for any finitely generated subgroup $H$ of $G$. Then,
$G \in {\scriptstyle{{\bf LH}}}\mathfrak{C}_{\beta}$. \hfill $\Box$

\vspace{0.1in}

\noindent
The definition of groups of type $\Phi$ by the second author in
\cite{T} is a particular case of a group class operation, applied
to the class $\mathfrak{F}$ of finite groups. More precisely, we
say that a group $G$ is of type $\Phi$ if the following conditions
are equivalent for a $\mathbb{Z}G$-module $M$:

(a) $\mbox{pd}_{\mathbb{Z}G}M < \infty$,

(b) $\mbox{pd}_{\mathbb{Z}H}\mbox{res}_H^GM < \infty$ for any
finite subgroup $H \subseteq G$.
\newline
Since the finitistic dimension of the integral group ring of a
finite group is equal to $1$, condition (b) above is equivalent
to either one of the following two conditions:

(b') $\mbox{pd}_{\mathbb{Z}H}\mbox{res}_H^GM \leq 1$ for any
finite subgroup $H \subseteq G$.

(b'') there is a non-negative integer $n=n(M)$, such that
$\mbox{pd}_{\mathbb{Z}H}\mbox{res}_H^GM \leq n$ for any
finite subgroup $H \subseteq G$.
\newline
More generally, let $k$ be a fixed commutative ring and consider
a group class $\mathfrak{C}$. Then, we let
$\Phi_{proj} \mathfrak{C}$ be the class consisting of
those groups $G$, over which the $kG$-modules of finite
projective dimension are precisely those $kG$-modules $M$
for which there is an integer $n=n(M) \geq 0$, such that
$\mbox{pd}_{kH}\mbox{res}_H^GM \leq n$ for any
$\mathfrak{C}$-subgroup $H \subseteq G$. It is clear that a
$kG$-module $M$ of finite projective dimension satisfies
the above condition with $n = \mbox{pd}_{kG}M$. (Compare
with \cite[Definition 5.6]{Bi1}.) The defining property of
a $\Phi_{proj} \mathfrak{C}$-group $G$ can be alternatively
described as follows:

\smallskip

\noindent
{\em If $M$ is a $kG$-module and
$\mbox{res}_H^GM \! \in {\tt Proj}(kH)$ for any
$\mathfrak{C}$-subgroup $H \subseteq G$, then
$\mbox{pd}_{kG}M < \infty$.}
\smallskip

\noindent
If $k = \mathbb{Z}$, then the groups of type $\Phi$ are precisely
the $\Phi_{proj}\mathfrak{F}$-groups. We also consider the variation
of the operation $\Phi_{proj}$ for flat modules and define for any
group class $\mathfrak{C}$ the class $\Phi_{flat} \mathfrak{C}$ to
be the class consisting of those groups $G$, over which the
$kG$-modules of finite flat dimension are precisely those
$kG$-modules $M$ that have uniformly finite flat dimension over $kH$
for any $\mathfrak{C}$-subgroup $H \subseteq G$ (with a uniform bound
that may depend on $M$). Alternatively, $G \in \Phi_{flat} \mathfrak{C}$
if and only if the following condition holds:

\smallskip

\noindent
{\em If $M$ is a $kG$-module and
$\mbox{res}_H^GM \in {\tt Flat}(kH)$ for any
$\mathfrak{C}$-subgroup $H \subseteq G$, then
$\mbox{fd}_{kG}M < \infty$.}
\smallskip

\noindent
Analogously, we may consider the operation $\Phi_{inj}$ and define
for any group class $\mathfrak{C}$ the class
$\Phi_{inj} \mathfrak{C}$ to be the class consisting of those
groups $G$, over which the $kG$-modules of finite injective
dimension are precisely those $kG$-modules $M$ that have
uniformly finite injective dimension over $kH$ for any
$\mathfrak{C}$-subgroup $H \subseteq G$ (with a uniform bound
that may depend on $M$). In the same way as before,
$G \in \Phi_{inj} \mathfrak{C}$ if and only if the following
condition holds:

\smallskip

\noindent
{\em If $M$ is a $kG$-module and
$\mbox{res}_H^GM \in {\tt Inj}(kH)$ for any
$\mathfrak{C}$-subgroup $H \subseteq G$, then
$\mbox{id}_{kG}M < \infty$.}

\smallskip

\begin{Lemma}
Let $k$ be a commutative ring. Then, the operations $\Phi_{proj}$,
$\Phi_{flat}$ and $\Phi_{inj}$ defined above are continuous.
\end{Lemma}
\vspace{-0.05in}
\noindent
{\em Proof.}
The proof of the continuity of the three operations is
analogous to that of Lemma 1.3 and is based on the fact
that the subgroups of a given group $G$ form a {\em set}.
\hfill $\Box$

\vspace{0.1in}

\noindent
We also note the following two properties of the operations
$\Phi_{proj}$, $\Phi_{flat}$ and $\Phi_{inj}$. The first of
these properties is an immediate consequence of Lambek's
flatness criterion \cite{L}: A module $M$ over a ring $R$
is flat if and only if its Pontryagin dual $DM$ is injective
(as a right $R$-module).

\begin{Lemma}
If $k$ is a commutative ring and $\mathfrak{C}$ is a group class,
then $\Phi_{inj}\mathfrak{C} \subseteq \Phi_{flat}\mathfrak{C}$.
\end{Lemma}
\vspace{-0.05in}
\noindent
{\em Proof.}
Let $G$ be a $\Phi_{inj}\mathfrak{C}$-group. In order to show
that $G \in \Phi_{flat}\mathfrak{C}$, consider a $kG$-module
$M$ and assume that $n$ is a non-negative integer, such that
$\mbox{fd}_{kH}\mbox{res}_H^GM \leq n$ for any
$\mathfrak{C}$-subgroup $H \subseteq G$. Then, for the dual
$kG$-module $DM$ we have $\mbox{res}_H^GDM = D\mbox{res}_H^GM$
and hence $\mbox{id}_{kH}\mbox{res}_H^GDM \leq n$ for any
$\mathfrak{C}$-subgroup $H \subseteq G$. Since
$G \in \Phi_{inj}\mathfrak{C}$, we conclude that
$\mbox{id}_{kG}DM < \infty$. It follows that
$\mbox{fd}_{kG}M < \infty$, as needed. \hfill $\Box$

\vspace{0.1in}

\noindent
Let $k$ be a commutative ring and consider a group $G$ and a
subgroup $H \subseteq G$. As shown by Gedrich and Gruenberg
in \cite{GG}, we have inequalities
$\mbox{spli} \, kH \leq \mbox{spli} \, kG$ and
$\mbox{silp} \, kH \leq \mbox{silp} \, kG$. In particular,
if $kG$ is Gorenstein regular, then $kH$ is Gorenstein
regular as well.

\begin{Proposition}
Let $k$ be a commutative ring and consider a group $G$, such
that $kG$ is Gorenstein regular. Then, the following conditions
are equivalent for a group class $\mathfrak{C}$:

(i) $G \in \Phi_{proj}\mathfrak{C}$,

(ii) $G \in \Phi_{flat}\mathfrak{C}$ and

(iii) $G \in \Phi_{inj}\mathfrak{C}$.
\end{Proposition}
\vspace{-0.05in}
\noindent
{\em Proof.}
For modules over a Gorenstein regular ring, the finiteness
of the projective dimension, the finiteness of the flat
dimension and the finiteness of the injective dimension are
conditions equivalent to each other. On the other hand, if
$kG$ is Gorenstein regular, then $kH$ is Gorenstein regular
for any subgroup $H \subseteq G$ (and, in particular, for any
$\mathfrak{C}$-subgroup $H \subseteq G$). Hence, the equivalence
of conditions (i), (ii) and (iii) is an immediate consequence
of the definitions. \hfill $\Box$

\section{Cokernels of acyclic complexes of projective/flat modules}

\noindent
In this Section, we fix a ring $R$ and consider, unless otherwise
specified, only left $R$-modules. Our goal is to relate the modules
that appear as cokernels of acyclic complexes of flat modules to
those modules that appear as cokernels of acyclic complexes of
projective modules.

Neeman \cite{N} proved that the embedding of the homotopy category
of projective modules into the homotopy category of flat modules
admits a right adjoint, whose kernel is the homotopy category of
pure-acyclic complexes of flat modules. He also showed that an
acyclic complex of flat modules is pure-acyclic if and only if
its cokernels are also flat modules. Therefore, for any complex
of flat modules $F$ there exists a complex of projective modules
$P$ and a chain map $f : P \longrightarrow F$, whose mapping cone
is pure-acyclic. We consider the induced short exact sequence of
chain complexes
\begin{equation}
 0 \longrightarrow F \longrightarrow \mbox{cone}(f)
   \longrightarrow P[1] \longrightarrow 0
\end{equation}
and note that, if $F$ is acyclic, then $P$ (and hence $P[1]$) is
also acyclic. Let $\mathcal{P}(R)$ (resp.\ $\mathcal{F}(R)$) be
the class consisting of those modules that appear as cokernels
of acyclic complexes of projective (resp.\ flat) modules. Then,
${\tt Proj}(R) \subseteq \mathcal{P}(R)$,
${\tt Flat}(R) \subseteq \mathcal{F}(R)$ and
$\mathcal{P}(R) \subseteq \mathcal{F}(R)$. We also consider the
class $\mathcal{I}(R)$ consisting of those modules that appear
as kernels of acyclic complexes of injective modules. Then,
${\tt Inj}(R) \subseteq \mathcal{I}(R)$ and
$D\mathcal{F}(R) \subseteq \mathcal{I}(R^{op})$. We note that
the classes $\mathcal{P}(R)$ and $\mathcal{F}(R)$ are closed
under direct sums, whereas $\mathcal{I}(R)$ is closed under
direct products.

\begin{Proposition}
For any $M \in \mathcal{F}(R)$ there exists a short exact sequence
of $R$-modules
\[ 0 \longrightarrow M \longrightarrow K \longrightarrow N
     \longrightarrow 0 , \]
where $K$ is flat and $N \in \mathcal{P}(R)$.
\end{Proposition}
\vspace{-0.05in}
\noindent
{\em Proof.}
We fix an acyclic complex of flat modules $F$ with $M=C_0(F)$ and
consider a short exact sequence of complexes as in (1) above. Then,
the complex $P[1]$ is acyclic and the induced short exact sequence
of $R$-modules
\[ 0 \longrightarrow C_0(F) \longrightarrow C_0(\mbox{cone}(f))
     \longrightarrow C_0(P[1]) \longrightarrow 0 \]
is of the required type. \hfill $\Box$

\vspace{0.1in}

\noindent
As an immediate consequence, we obtain the equivalence between
the triviality of all acyclic complexes of projective modules
and the triviality of all acyclic complexes of flat modules.

\begin{Corollary}
The following conditions are equivalent for a ring $R$:

(i) $\mathcal{P}(R) = {\tt Proj}(R)$, i.e.\ any acyclic complex
of projective modules is contractible.

(ii) $\mathcal{F}(R) = {\tt Flat}(R)$, i.e.\ any acyclic complex
of flat modules is pure-acyclic.
\newline
These conditions are satisfied if
$\mathcal{I}(R^{op}) = {\tt Inj}(R^{op})$, i.e.\ if any acyclic
complex of injective right modules is contractible.
\end{Corollary}
\vspace{-0.05in}
\noindent
{\em Proof.}
(i)$\rightarrow$(ii): Let $M$ be an $\mathcal{F}(R)$-module and
consider a short exact sequence as in Proposition 2.1. Assumption
(i) implies that the module $N$ therein is projective and hence
the short exact sequence splits. It follows that $M$ is a direct
summand of a flat module, i.e.\ $M$ is flat.

(ii)$\rightarrow$(i): Let $P$ be an acyclic complex of projective
modules. Since $P$ is also an acyclic complex of flat modules,
assumption (ii) implies that the cokernels of $P$ are all flat.
Then, the contractibility of $P$ follows from \cite{BG} (see also
\cite{N}).

Finally, assume that $\mathcal{I}(R^{op}) = {\tt Inj}(R^{op})$.
Then, we have an inclusion
$D\mathcal{F}(R) \subseteq {\tt Inj}(R^{op})$ and hence
$\mathcal{F}(R) \subseteq D^{-1}{\tt Inj}(R^{op}) = {\tt Flat}(R)$,
where the latter equality is precisely Lambek's flatness criterion
\cite{L}. \hfill $\Box$

\vspace{0.1in}

\noindent
The following result is an analogue of Corollary 2.2 in Gorenstein
homological algebra. We note that, in view of the very definition
of the Gorenstein modules, there are inclusions
${\tt PGF}(R) \subseteq \mathcal{P}(R)$,
${\tt GFlat}(R) \subseteq \mathcal{F}(R)$ and
${\tt GInj}(R) \subseteq \mathcal{I}(R)$.

\begin{Corollary}
The following conditions are equivalent for a ring $R$:

(i) $\mathcal{P}(R) = {\tt PGF}(R)$,

(ii) $\mathcal{F}(R) = {\tt GFlat}(R)$.
\newline
These conditions are satisfied if
$\mathcal{I}(R^{op}) = {\tt GInj}(R^{op})$.
\end{Corollary}
\vspace{-0.05in}
\noindent
{\em Proof.}
(i)$\rightarrow$(ii): Let $M$ be an $\mathcal{F}(R)$-module
and consider a short exact sequence as in Proposition 2.1.
Then, assumption (i) implies that the module $N$ therein is
(projectively coresolved) Gorenstein flat. Since flat modules
are also Gorenstein flat, the closure of ${\tt GFlat}(R)$ under
kernels of epimorphisms implies that $M$ is Gorenstein flat.

(ii)$\rightarrow$(i): Let $P$ be an acyclic complex of projective
modules. Since $P$ is, of course, an acyclic complex of flat modules,
assumption (ii) implies that the cokernels of $P$ are Gorenstein
flat. Then, the functor $\mbox{Tor}_R^1(I,\_\!\_)$ vanishes on these
cokernels for any injective right module $I$. It follows that the
complex of abelian groups $I \otimes_{kG}P$ is acyclic for any
injective right module $I$, so that the cokernels of $P$ are
projectively coresolved Gorenstein flat.

Finally, we assume that $\mathcal{I}(R^{op}) = {\tt GInj}(R^{op})$
and show that condition (ii) is satisfied. If $F$ is an acyclic
complex of flat modules, then $DF$ is an acyclic complex of injective
right modules and hence (in view of our assumption) it is totally
acyclic. The complex of abelian groups
$\mbox{Hom}_{R^{op}}(I,DF) = D(I \otimes_R F)$ is therefore acyclic
for any injective right module $I$. It follows that $I \otimes_R F$
is an acyclic complex of abelian groups for any injective right
module $I$, so that $F$ is a totally acyclic complex of flat
modules, as needed. \hfill $\Box$

\vspace{0.1in}

\noindent
We say that a module $M$ has finite sfp-injective dimension if
$M$ admits an injective resolution with some cosyzygy module
strongly fp-injective, in the sense of \cite{LGO}. As shown in
\cite{EK, ET1}, all modules of finite sfp-injective dimension
are contained in $\mathcal{P}(R)^{\perp}$ and all cotorsion
modules of finite sfp-injective dimension are contained in
$\mathcal{F}(R)^{\perp}$. In general, the relation between
these two right orthogonal classes is described below.

\begin{Corollary}
$\mathcal{F}(R)^{\perp} = {\tt Cotor}(R) \cap \mathcal{P}(R)^{\perp}$.
\end{Corollary}
\vspace{-0.05in}
\noindent
{\em Proof.}
Since $\mathcal{F}(R)$ contains all flat and all $\mathcal{P}(R)$-modules,
it is clear that $\mathcal{F}(R)^{\perp}$ is contained in both
${\tt Flat}(R)^{\perp} = {\tt Cotor}(R)$ and $\mathcal{P}(R)^{\perp}$.
Conversely, assume that $L$ is a cotorsion module contained in
$\mathcal{P}(R)^{\perp}$. We shall prove that
$L \in \mathcal{F}(R)^{\perp}$, i.e.\ that the functor
$\mbox{Ext}_R^1(\_\!\_,L)$ vanishes on $\mathcal{F}(R)$. To that end,
we fix a module $M \in \mathcal{F}(R)$ and consider a short exact
sequence as in Proposition 2.1. Since any $\mathcal{P}(R)$-module admits
a projective resolution with cokernels in $\mathcal{P}(R)$, it is easily
seen that the functors $\mbox{Ext}_R^i(\_\!\_,L)$ vanish on
$\mathcal{P}(R)$ for all $i \geq 1$. In particular, $\mbox{Ext}_R^2(N,L)=0$.
Since $\mbox{Ext}_R^1(K,L)=0$, we conclude that the abelian group
$\mbox{Ext}_R^1(M,L)$ is trivial, as needed. \hfill $\Box$

\vspace{0.1in}

\noindent
{\bf Remarks 2.5.}
(i) In view of Corollary 2.3, all acyclic complexes of projective
modules remain acyclic after applying the functor $I \otimes_R \_\!\_$
for all injective right modules $I$ if and only if the same holds for
all acyclic complexes of flat modules. These conditions are satisfied
if all acyclic complexes of injective right modules remain acyclic
after applying the functor $\mbox{Hom}_{R^{op}}(I,\_\!\_)$ for any
injective right module $I$. More generally, let $\mathcal{A}$ be any
class of right modules. Then, the existence of a short exact sequence
of complexes as in (1) for any (acyclic) complex of flat modules implies
that all acyclic complexes of projective modules remain acyclic after
applying the functor $A \otimes_R \_\!\_$ for any right module
$A \in \mathcal{A}$ if and only if the same holds for all acyclic
complexes of flat modules. These conditions are satisfied if all
acyclic complexes of injective right modules remain acyclic after
applying the functor $\mbox{Hom}_{R^{op}}(A,\_\!\_)$ for any right
module $A \in \mathcal{A}$.

(ii) Let $M \in \mathcal{F}(R)$ and consider a short exact sequence
of modules as in Proposition 2.1. We also consider a short exact
sequence
\[ 0 \longrightarrow K' \longrightarrow Q \longrightarrow K
     \longrightarrow 0 , \]
where $Q$ is projective. The pullback of this short exact sequence
along the monomorphism $M \longrightarrow K$ provides us with another
short exact sequence
\[ 0 \longrightarrow K' \longrightarrow N' \longrightarrow M
     \longrightarrow 0 , \]
where $K'$ is flat and $N' \in \mathcal{P}(R)$.

(iii) Proposition 2.1 implies that
$\mathcal{F}(R) \cap \mathcal{P}(R)^{\perp} \subseteq
 {\tt Flat}(R)$.
Assume that any flat module has finite sfp-injective
dimension.\footnote{This is the case if, for example,
any flat module has finite injective dimension, i.e.\
if $\mbox{silp} \, R$ is finite; cf.\
\cite[Proposition 2.1]{ET0}.} Then, all flat modules
are contained in $\mathcal{P}(R)^{\perp}$ and hence
$\mathcal{F}(R) \cap \mathcal{P}(R)^{\perp} = {\tt Flat}(R)$.

\section{Acyclic complexes of projective/flat modules over group algebras}

\noindent
In this Section, we consider those groups $G$, for which the cokernels
of any acyclic complex of projective (resp.\ flat) modules over the
group algebra $kG$ are Gorenstein projective (resp.\ Gorenstein flat).
We work over a fixed commutative ring $k$ and note that for any group
$G$ there are inclusions
\begin{equation}
 {\tt Proj}(kG) \subseteq  {\tt PGF}(kG) \subseteq
 {\tt GProj}(kG) \subseteq \mathcal{P}(kG)
 \;\; \mbox{and} \;\;
 {\tt Flat}(kG) \subseteq  {\tt GFlat}(kG) \subseteq
 \mathcal{F}(kG) .
\end{equation}

\begin{Lemma}
The following conditions are equivalent for a group $G$:

(i) ${\tt PGF}(kG) = {\tt GProj}(kG) = \mathcal{P}(kG)$,

(ii) $\mathcal{P}(kG) \subseteq {\tt PGF}(kG)$,

(iii) ${\tt GFlat}(kG) = \mathcal{F}(kG)$ and

(iv) $\mathcal{F}(kG) \subseteq {\tt GFlat}(kG)$.
\newline
In addition, the class $\mathfrak{X} = \mathfrak{X}(k)$ consisting
of those group $G$ that satisfy these conditions, is subgroup-closed.
\end{Lemma}
\vspace{-0.05in}
\noindent
{\em Proof.}
In view of the inclusions (2), it is clear that (i)$\leftrightarrow$(ii)
and (iii)$\leftrightarrow$(iv). Since (i) is satisfied if and only if
${\tt PGF}(kG) = \mathcal{P}(kG)$, the equivalence (i)$\leftrightarrow$(iii)
follows from Corollary 2.3.

It remains to prove that the class $\mathfrak{X}$ is subgroup-closed.
To that end, consider an $\mathfrak{X}$-group $G$ and let $H \subseteq G$
be a subgroup. We have to show that $H$ is an $\mathfrak{X}$-group as
well, i.e.\ that $\mathcal{P}(kH) \subseteq {\tt PGF}(kH)$. Let $P$ be
an acyclic complex of projective $kH$-modules and consider an injective
$kH$-module $I$. We aim at proving that the complex of abelian groups
$I \otimes_{kH}P$ is acyclic. In fact, since $I$ is a direct summand of
$\mbox{res}_H^G\mbox{coind}_H^GI$, the acyclicity of $I \otimes_{kH}P$
will follow if we show that the complex
$\mbox{res}_H^G\mbox{coind}_H^GI \otimes_{kH} P$ is acyclic. Now,
the associativity of the tensor product implies the existence of
an isomorphism of complexes
\[ \mbox{res}_H^G\mbox{coind}_H^GI \otimes_{kH} P \simeq
   \mbox{coind}_H^GI \otimes_{kG} \mbox{ind}_H^GP . \]
Since $\mbox{ind}_H^GP$ is an acyclic complex of projective
$kG$-modules and $\mbox{coind}_H^GI$ is an injective $kG$-module,
our hypothesis on $G$, namely that
$\mathcal{P}(kG) \subseteq {\tt PGF}(kG)$, implies that
$\mbox{res}_H^G\mbox{coind}_H^GI \otimes_{kH} P$ is an acyclic
complex, as needed. \hfill $\Box$

\vspace{0.1in}

\noindent
Groups in class $\mathfrak{X}$ have (essentially by their definition)
a smooth behaviour with respect to Gorenstein homological algebra notions.

\begin{Proposition}
If $G$ is an $\mathfrak{X}$-group, then:

(i) The cotorsion pair
$\left( {\tt GProj}(kG) , {\tt GProj}(kG)^{\perp} \right)$ is
complete.

(ii) ${\tt GProj}(kG) \subseteq {\tt GFlat}(kG)$.

(iii) Any acyclic complex of projective $kG$-modules is totally
acyclic, i.e.\ it remains acyclic after applying the functor
$\mbox{Hom}_{kG}(\_\!\_,P)$ for any projective $kG$-module $P$.

(iv) Any acyclic complex of flat $kG$-modules is totally acyclic,
i.e.\ it remains acyclic after applying the functor
$I \otimes_{kG} \_\!\_$ for any injective $kG$-module $I$.
\end{Proposition}
\vspace{-0.05in}
\noindent
{\em Proof.}
Since $G \in \mathfrak{X}$, we know that
${\tt PGF}(kG) = {\tt GProj}(kG)$, Then, (i) follows since the
cotorsion pair $\left( {\tt PGF}(kG),{\tt PGF}(kG)^{\perp} \right)$
is complete, whereas (ii) follows since
${\tt PGF}(kG) \subseteq {\tt GFlat}(kG)$. We also note that
assertion (iii) is a restatement of the equality
${\tt GProj}(kG) = \mathcal{P}(kG)$ and (iv) is a restatement
of the equality ${\tt GFlat}(kG) = \mathcal{F}(kG)$. \hfill $\Box$

\vspace{0.1in}

\noindent
In order to provide a supply of examples of $\mathfrak{X}$-groups,
we shall use the following result.

\begin{Theorem}
The class $\mathfrak{X}$ is (i) ${\scriptstyle{{\bf LH}}}$-closed,
(ii) $\Phi_{proj}$-closed and (iii) $\Phi_{flat}$-closed.
\end{Theorem}
\vspace{-0.05in}
\noindent
{\em Proof.}
(i) We have to show that
${\scriptstyle{{\bf LH}}}\mathfrak{X} = \mathfrak{X}$, i.e.\
that
${\scriptstyle{{\bf LH}}}\mathfrak{X} \subseteq \mathfrak{X}$.
Let $G$ be a group contained in
${\scriptstyle{{\bf LH}}}\mathfrak{X}$. In order to show that
$G$ is an $\mathfrak{X}$-group, we shall prove that
$\mathcal{P}(kG) \subseteq {\tt PGF}(kG)$.

First of all, we show that for all $M \in \mathcal{P}(kG)$ and
all ${\scriptstyle{\bf H}}\mathfrak{X}$-subgroups $H \subseteq G$
the $kH$-module $\mbox{res}_H^GM$ is contained in ${\tt PGF}(kH)$.
To that end, we use transfinite induction on the ordinal $\alpha$,
which is such that
$H \in {\scriptstyle{{\bf H}}}_{\alpha}\mathfrak{X}$. If
$\alpha = 0$, then $H \in \mathfrak{X}$. Since
$\mbox{res}_H^GM \in \mathcal{P}(kH)$, it follows that
$\mbox{res}_H^GM \in {\tt PGF}(kH)$. We now assume that $\alpha >0$
and the result is true for all $kG$-modules $M \in \mathcal{P}(kG)$
and all ${\scriptstyle{{\bf H}}}_{\beta}\mathfrak{X}$-subgroups
$H \subseteq G$ for all $\beta < \alpha$. Consider
a $kG$-module $M \in \mathcal{P}(kG)$ and let $H$ be an
${\scriptstyle{{\bf H}}}_{\alpha}\mathfrak{X}$-subgroup
of $G$. Then, there exists an exact sequence of $kH$-modules
\[ 0 \longrightarrow M_r \longrightarrow \cdots
     \longrightarrow M_1 \longrightarrow M_0
     \longrightarrow \mbox{res}_H^GM \longrightarrow 0 , \]
where each $M_i$ is a direct sum of $kH$-modules of the form
$\mbox{ind}_F^H \mbox{res}_F^GM$, for an
${\scriptstyle{{\bf H}}}_{\beta}\mathfrak{X}$-subgroup
$F \subseteq H$ for some $\beta < \alpha$. We note that $r$ does not depend
on $M$; it is simply the dimension of the $H$-CW-complex witnessing
that $H \in {\scriptstyle{{\bf H}}}_{\alpha}\mathfrak{X}$. Since
induction from $F$ to $H$ maps ${\tt PGF}(kF)$ into ${\tt PGF}(kH)$
for any subgroup $F \subseteq H$, we conclude that $M_i \in {\tt PGF}(kH)$
for all $i=0,1, \ldots ,r$ and hence the $kH$-module $\mbox{res}_H^GM$
has PGF-dimension $\leq r$. Since this holds for any
$M \in \mathcal{P}(kG)$, a simple argument shows that the $kH$-module
$\mbox{res}_H^GM$ is actually contained in ${\tt PGF}(kH)$. Indeed,
if $P$ is an acyclic complex of projective $kG$-modules witnessing
that $M \in \mathcal{P}(kG)$, then $M$ is the $r$-th syzygy of another
$kG$-module $M' \in \mathcal{P}(kG)$. It follows that $\mbox{res}_H^GM$
is the $r$-th syzygy of $\mbox{res}_H^GM'$ and hence we conclude that
$\mbox{res}_H^GM \in {\tt PGF}(kH)$.

We shall now prove that $\mbox{res}_H^GM \in {\tt PGF}(kH)$ for all
$M \in \mathcal{P}(kG)$ and all
${\scriptstyle{\bf LH}}\mathfrak{X}$-subgroups $H \subseteq G$. (In
particular, letting $H=G$, it will follow that
$\mathcal{P}(kG) \subseteq {\tt PGF}(kG)$.) We proceed by induction
on the cardinality $\kappa$ of $H$. If $\kappa \leq \aleph_0$, then the
${\scriptstyle{{\bf LH}}}\mathfrak{X}$-group $H$ is actually contained
in ${\scriptstyle{{\bf H}}}\mathfrak{X}$ and we are done by the previous
discussion. If $\kappa$ is uncountable, we may express $H$ as a continuous
ascending union of subgroups $(H_{\alpha})_{\alpha < \lambda}$ of cardinality
$< \kappa$, which is indexed by a suitable ordinal $\lambda$. The class
$\mathfrak{X}$ is subgroup-closed (cf.\ Lemma 3.1) and hence this is also
the case for the class ${\scriptstyle{\bf LH}}\mathfrak{X}$; in particular,
$H_{\alpha} \in {\scriptstyle{\bf LH}}\mathfrak{X}$ for all $\alpha$.
Invoking our induction hypothesis, it follows that $\mbox{res}_{H_{\alpha}}^GM$
is contained in ${\tt PGF}(kH_{\alpha})$ for all $\alpha$. Then,
we may conclude as in \cite[Lemma 5.6]{Ben} that the $kH$-module
$\mbox{res}_H^GM$ has PGF-dimension $\leq 1$. Since this is the case
for any $M \in \mathcal{P}(kG)$, we can show as above that the
$kH$-module $\mbox{res}_H^GM$ is contained in ${\tt PGF}(kH)$.

(ii) We have to show that
$\Phi_{proj}\mathfrak{X} = \mathfrak{X}$, i.e.\ that
$\Phi_{proj}\mathfrak{X} \subseteq \mathfrak{X}$. Let $G$ be a group
contained in $\Phi_{proj}\mathfrak{X}$. In order to show that $G$ is
an $\mathfrak{X}$-group, we shall prove that
$\mathcal{P}(kG) \subseteq {\tt PGF}(kG)$. To that end, we consider a
$kG$-module $M \in \mathcal{P}(kG)$ and let
\begin{equation}
 0 \longrightarrow N \longrightarrow P
   \longrightarrow M \longrightarrow 0
\end{equation}
be a short exact sequence of $kG$-modules, where
$N \in {\tt PGF}(kG)^{\perp}$ and $P \in {\tt PGF}(kG)$. We fix an
$\mathfrak{X}$-subgroup $H \subseteq G$. Since $M,P \in \mathcal{P}(kG)$,
it follows that $\mbox{res}_H^GM, \mbox{res}_H^GP \in \mathcal{P}(kH)$
and hence $\mbox{res}_H^GM, \mbox{res}_H^GP \in {\tt PGF}(kH)$. The class
${\tt PGF}(kH)$ is closed under kernels of epimorphisms and hence the
$kH$-module $\mbox{res}_H^GN$ is also contained in ${\tt PGF}(kH)$.
On the other hand, induction from $H$ to $G$ maps ${\tt PGF}(kH)$ into
${\tt PGF}(kG)$ and hence the induction-restriction isomorphism
\[ \mbox{Ext}^1_{kG}(\mbox{ind}_H^G \_\!\_ ,N) \simeq
   \mbox{Ext}^1_{kH}(\_\!\_ , \mbox{res}_H^GN) \]
shows that $\mbox{res}_H^GN \in {\tt PGF}(kH)^{\perp}$. We conclude
that $\mbox{res}_H^GN \in {\tt PGF}(kH) \cap {\tt PGF}(kH)^{\perp}$
and hence the $kH$-module $\mbox{res}_H^GN$ is projective. Since this
is the case for any $\mathfrak{X}$-subgroup $H \subseteq G$, our
assumption that $G$ is a $\Phi_{proj} \mathfrak{X}$-group implies
that $\mbox{pd}_{kG}N < \infty$; in particular, the $kG$-module $N$
has finite PGF-dimension. Then, it follows from the short exact
sequence (3) that $M$ has finite PGF-dimension as well. Having proved
that any $kG$-module $M \in \mathcal{P}(kG)$ has finite PGF-dimension,
the closure of $\mathcal{P}(kG)$ under direct sums implies that there
is an upper bound on the PGF-dimensions of $\mathcal{P}(kG)$-modules.
Then, arguing as in (i), we can show that any $\mathcal{P}(kG)$-module
is actually contained in ${\tt PGF}(kG)$.

(iii) We have to show that $\Phi_{flat}\mathfrak{X} = \mathfrak{X}$,
i.e.\ that $\Phi_{flat}\mathfrak{X} \subseteq \mathfrak{X}$. Let $G$
be a group contained in $\Phi_{flat}\mathfrak{X}$. In order to show
that $G$ is an $\mathfrak{X}$-group, we shall prove that
$\mathcal{F}(kG) \subseteq {\tt GFlat}(kG)$; cf.\ Lemma 3.1. The
argument is similar to that used in (ii) above. We consider a
$kG$-module $M \in \mathcal{F}(kG)$ and let
\begin{equation}
 0 \longrightarrow N \longrightarrow F
   \longrightarrow M \longrightarrow 0
\end{equation}
be a short exact sequence of $kG$-modules, where
$N \in {\tt GFlat}(kG)^{\perp}$ and $F \in {\tt GFlat}(kG)$. We fix an
$\mathfrak{X}$-subgroup $H \subseteq G$. Since $M,F \in \mathcal{F}(kG)$,
it follows that $\mbox{res}_H^GM, \mbox{res}_H^GF \in \mathcal{F}(kH)$
and hence $\mbox{res}_H^GM, \mbox{res}_H^GF \in {\tt GFlat}(kH)$. Then,
the closure of the class ${\tt GFlat}(kH)$ under kernels of epimorphisms
implies that the $kH$-module $\mbox{res}_H^GN$ is Gorenstein flat as well.
On the other hand, induction from $H$ to $G$ maps ${\tt GFlat}(kH)$ into
${\tt GFlat}(kG)$ and hence the induction-restriction isomorphism
\[ \mbox{Ext}^1_{kG}(\mbox{ind}_H^G \_\!\_ ,N) \simeq
   \mbox{Ext}^1_{kH}(\_\!\_ , \mbox{res}_H^GN) \]
shows that $\mbox{res}_H^GN \in {\tt GFlat}(kH)^{\perp}$. We conclude
that $\mbox{res}_H^GN \in {\tt GFlat}(kH) \cap {\tt GFlat}(kH)^{\perp}$,
i.e.\ the $kH$-module $\mbox{res}_H^GN$ is flat-cotorsion. Since this is
the case for any $\mathfrak{X}$-subgroup $H \subseteq G$, our assumption
that $G \in \Phi_{flat} \mathfrak{X}$ implies that
$\mbox{fd}_{kG}N < \infty$; in particular, $N$ has finite Gorenstein flat
dimension. Then, the short exact sequence (4) shows that $M$ has finite
Gorenstein flat dimension as well. Having proved that any $kG$-module
$M \in \mathcal{F}(kG)$ has finite Gorenstein flat dimension, the closure
of $\mathcal{F}(kG)$ under direct sums implies that there is an upper
bound on the Gorenstein flat dimensions of all $\mathcal{F}(kG)$-modules.
Then, arguing as in (i), we can show that any $\mathcal{F}(kG)$-module is
actually Gorenstein flat. \hfill $\Box$

\vspace{0.1in}

\noindent
As we have alluded to in the statement of Lemma 3.1, the class
$\mathfrak{X}$ depends upon the coefficient ring $k$. In order
for the trivial group to be contained in $\mathfrak{X}$, it is
necessary that $\mathcal{P}(k)={\tt PGF}(k)$ or, equivalently,
that $\mathcal{F}(k) = {\tt GFlat}(k)$. The latter conditions
are satisfied if, for example, $k$ is weakly Gorenstein regular,
i.e.\ if $\mbox{sfli} \, k$ is finite. (Indeed, if $M$ is a
$k$-module of finite flat dimension, then the functor
$M \otimes_k \_\!\_$ preserves the acyclicity of any acyclic
complex of flat $k$-modules.)

\noindent
\begin{Proposition}
Assume that $\mathcal{P}(k) = {\tt PGF}(k)$ or, equivalently, that
$\mathcal{F}(k) = {\tt GFlat}(k)$. Then, all finite groups are
contained in $\mathfrak{X}$.
\end{Proposition}
\vspace{-0.05in}
\noindent
{\em Proof.}
We have to show that $\mathcal{P}(kG) \subseteq {\tt PGF}(kG)$, in the
special case where $G$ is a finite group. In other words, we have to
show that any acyclic complex $P$ of projective $kG$-modules remains
acyclic upon applying the functor $I \otimes_{kG} \_\!\_$ for any
injective $kG$-module $I$. Since any injective $kG$-module is a
direct summand of a coinduced $kG$-module of the form
$\mbox{coind}_1^GJ \simeq \mbox{ind}_1^GJ$, where $J$ is an injective
$k$-module, we are reduced to showing that the complex
$\mbox{ind}_1^GJ \otimes_{kG} P = J \otimes_k \mbox{res}_1^GP$ is
acyclic. The latter claim follows from our assumption that
$\mathcal{P}(k) = {\tt PGF}(k)$. \hfill $\Box$

\vspace{0.1in}

\noindent
Let $\mathfrak{X}_{fin}$ be the smallest group class which contains
all finite groups and is ${\scriptstyle{{\bf LH}}}$-closed,
$\Phi_{proj}$-closed and $\Phi_{flat}$-closed. We note that groups
in class $\mathfrak{X}_{fin}$ admit a hierarchical description with
base the class of finite groups, as explained in $\S $1.IV.

\begin{Corollary}
Assume that $\mathcal{P}(k) = {\tt PGF}(k)$ or, equivalently,
that $\mathcal{F}(k) = {\tt GFlat}(k)$ and let $G$ be an
$\mathfrak{X}_{fin}$-group. Then:

(i) The cotorsion pair
$\left( {\tt GProj}(kG) , {\tt GProj}(kG)^{\perp} \right)$ is
complete.

(ii) ${\tt GProj}(kG) \subseteq {\tt GFlat}(kG)$.

(iii) Any acyclic complex of projective $kG$-modules is totally
acyclic, i.e.\ it remains acyclic after applying the functor
$\mbox{Hom}_{kG}(\_\!\_,P)$ for any projective $kG$-module $P$.

(iv) Any acyclic complex of flat $kG$-modules is totally acyclic,
i.e.\ it remains acyclic after applying the functor
$I \otimes_{kG} \_\!\_$ for any injective $kG$-module $I$.
\end{Corollary}
\vspace{-0.05in}
\noindent
{\em Proof.}
Since $\mathfrak{X}$ contains the class of finite groups (cf.\
Proposition 3.4), Theorem 3.3 implies that
$\mathfrak{X}_{fin} \subseteq \mathfrak{X}$. The result is
therefore an immediate consequence of Proposition 3.2. \hfill $\Box$

\vspace{0.1in}

\noindent
{\bf Remarks 3.6.}
(i) Let $k$ be a commutative ring. As we noted before,
the equivalent conditions $\mathcal{P}(k) = {\tt PGF}(k)$
and $\mathcal{F}(k) = {\tt GFlat}(k)$ are satisfied if
$k$ is weakly Gorenstein regular. Hence, the pairs $(k,G)$
considered in Corollary 3.5 include those considered in
\cite{Bi2}, where it was assumed that $k$ has finite global
dimension and $G$ is either an
${\scriptstyle{{\bf LH}}}\mathfrak{F}$-group or a group of
type $\Phi$.

(ii) Assume that $k = \mathbb{Z}$ is the ring of integers.
It would be of interest to know of a group $G$ which is
contained in $\mathfrak{X}_{fin}$, but is not an
${\scriptstyle{{\bf LH}}}\mathfrak{F}$-group nor a group
of type $\Phi$.
\addtocounter{Lemma}{1}

\section{Acyclic complexes of injective modules over group algebras}

\noindent
In this section, we consider groups $G$, for which the kernels of
any acyclic complex of injective modules over $kG$ are Gorenstein
injective. Having fixed the commutative ring $k$, we note that for
any group $G$ there are inclusions
\[ {\tt Inj}(kG) \subseteq  {\tt GInj}(kG) \subseteq
   \mathcal{I}(kG) . \]
We consider the group class $\mathfrak{Y} = \mathfrak{Y}(k)$,
which is defined by letting
\[ \mathfrak{Y} = \{ G : \mathcal{I}(kG) \subseteq {\tt GInj}(kG) \} =
   \{ G : {\tt GInj}(kG) = \mathcal{I}(kG) \} . \]
Then, a group $G$ is contained in $\mathfrak{Y}$ if and only if
any acyclic complex of injective $kG$-modules is totally acyclic,
i.e.\ if and only if any acyclic complex of injective $kG$-modules
remains acyclic after applying the functor $\mbox{Hom}_{kG}(I,\_\!\_)$
for any injective $kG$-module $I$. In view of Corollary 2.3,
$\mathfrak{Y}$ is a subclass of the class $\mathfrak{X}$ introduced
in the previous section.

\begin{Proposition}
The class $\mathfrak{Y}$ is subgroup-closed and for any
$\mathfrak{Y}$-group $G$ the following hold:

(i) The cotorsion pair
$\left( {\tt GProj}(kG) , {\tt GProj}(kG)^{\perp} \right)$ is
complete.

(ii) ${\tt GProj}(kG) \subseteq {\tt GFlat}(kG)$.

(iii) Any acyclic complex of projective $kG$-modules is totally
acyclic, i.e.\ it remains acyclic after applying the functor
$\mbox{Hom}_{kG}(\_\!\_,P)$ for any projective $kG$-module $P$.

(iv) Any acyclic complex of flat $kG$-modules is totally acyclic,
i.e.\ it remains acyclic after applying the functor
$I \otimes_{kG} \_\!\_$ for any injective $kG$-module $I$.

(v) Any acyclic complex of injective $kG$-modules is totally
acyclic, i.e.\ it remains acyclic after applying the functor
$\mbox{Hom}_{kG}(I,\_\!\_)$ for any injective $kG$-module $I$.
\end{Proposition}
\vspace{-0.05in}
\noindent
{\em Proof.}
As we noted above, $\mathfrak{Y}$ is a subclass of $\mathfrak{X}$,
so that assertions (i)-(iv) follow from Proposition 3.2. On the
other hand, assertion (v) is precisely the definition of a
$\mathfrak{Y}$-group.

It remains to show that $\mathfrak{Y}$ is subgroup-closed. To that
end, let $G$ be a $\mathfrak{Y}$-group and $H \subseteq G$ be a
subgroup. In order to show that $H$ is a $\mathfrak{Y}$-group as
well, we consider an acyclic complex of injective $kH$-modules $J$,
an injective $kH$-module $I$ and aim at proving that the complex
of abelian groups $\mbox{Hom}_{kH}(I,J)$ is acyclic. Since $I$ is
a direct summand of $\mbox{res}_H^G\mbox{coind}_H^GI$, it suffices
to show that the complex
$\mbox{Hom}_{kH}(\mbox{res}_H^G\mbox{coind}_H^GI,J) \simeq
 \mbox{Hom}_{kG}(\mbox{coind}_H^GI,\mbox{coind}_H^GJ)$
is acyclic. The latter claim follows from our hypothesis that
$G \in \mathfrak{Y}$, since $\mbox{coind}_H^GJ$ is an acyclic
complex of injective $kG$-modules and $\mbox{coind}_H^GI$ is an
injective $kG$-module. \hfill $\Box$

\vspace{0.1in}

\noindent
In order to provide a supply of examples of $\mathfrak{Y}$-groups,
we shall use the following result, whose proof is similar to that
of Theorem 3.3.

\begin{Theorem}
The class $\mathfrak{Y}$ is (i) ${\scriptstyle{{\bf LH}}}$-closed
and (ii) $\Phi_{inj}$-closed.
\end{Theorem}
\vspace{-0.05in}
\noindent
{\em Proof.}
(i) We have to show that
${\scriptstyle{{\bf LH}}}\mathfrak{Y} = \mathfrak{Y}$,
i.e.\ that
${\scriptstyle{{\bf LH}}}\mathfrak{Y} \subseteq \mathfrak{Y}$.
To that end, we consider an
${\scriptstyle{{\bf LH}}}\mathfrak{Y}$-group $G$ and prove that
$G \in \mathfrak{Y}$, by showing that
$\mathcal{I}(kG) \subseteq {\tt GInj}(kG)$.

First of all, we show that for all $M \in \mathcal{I}(kG)$ and
all ${\scriptstyle{\bf H}}\mathfrak{Y}$-subgroups $H \subseteq G$
the $kH$-module $\mbox{res}_H^GM$ is Gorenstein injective. To
that end, we use transfinite induction on the ordinal number
$\alpha$, which is such that
$H \in {\scriptstyle{{\bf H}}}_{\alpha}\mathfrak{Y}$. If
$\alpha = 0$, then $H \in \mathfrak{Y}$. Since
$\mbox{res}_H^GM \in \mathcal{I}(kH)$, it follows that
$\mbox{res}_H^GM \in {\tt GInj}(kH)$. We now assume that
$\alpha >0$ and the result is true for all $kG$-modules
$M \in \mathcal{I}(kG)$ and all
${\scriptstyle{{\bf H}}}_{\beta}\mathfrak{Y}$-subgroups
$H \subseteq G$ for all $\beta < \alpha$. We consider
a $kG$-module $M \in \mathcal{I}(kG)$ and let $H$ be an
${\scriptstyle{{\bf H}}}_{\alpha}\mathfrak{Y}$-subgroup
of $G$. Then, there exists an exact sequence of $kH$-modules
\[ 0 \longrightarrow \mbox{res}_H^GM \longrightarrow M^0
     \longrightarrow M^1 \longrightarrow \cdots
     \longrightarrow M^r \longrightarrow 0 , \]
where each $M^i$ is a direct product of $kH$-modules of the form
$\mbox{coind}_F^H \mbox{res}_F^GM$, for a suitable
${\scriptstyle{{\bf H}}}_{\beta}\mathfrak{Y}$-subgroup
$F \subseteq H$ for some $\beta < \alpha$. We note that $r$ does
not depend on the $kG$-module $M$; it is simply the dimension of
the $H$-CW-complex witnessing that
$H \in {\scriptstyle{{\bf H}}}_{\alpha}\mathfrak{Y}$. Since
coinduction from $F$ to $H$ maps the class ${\tt GInj}(kF)$ into
${\tt GInj}(kH)$ for any subgroup $F \subseteq H$, we conclude
that $M^i \in {\tt GInj}(kH)$ for all $i=0,1, \ldots ,r$ and
hence the $kH$-module $\mbox{res}_H^GM$ has Gorenstein injective
dimension $\leq r$. This is the case for any $M \in \mathcal{I}(kG)$.
A simple argument shows that the $kH$-module $\mbox{res}_H^GM$ is
actually Gorenstein injective. Indeed, if $I$ is an acyclic complex
of injective $kG$-modules witnessing that $M \in \mathcal{I}(kG)$,
then $M$ is the $r$-th cosyzygy of another $kG$-module
$M' \in \mathcal{I}(kG)$. It follows that $\mbox{res}_H^GM$ is the
$r$-th cosyzygy of $\mbox{res}_H^GM'$ and hence
$\mbox{res}_H^GM \in {\tt GInj}(kH)$.

We shall now prove that $\mbox{res}_H^GM \in {\tt GInj}(kH)$ for all
$M \in \mathcal{I}(kG)$ and all
${\scriptstyle{\bf LH}}\mathfrak{Y}$-subgroups $H \subseteq G$. (In
particular, letting $H=G$, it will follow that
$\mathcal{I}(kG) \subseteq {\tt GInj}(kG)$.) We proceed by induction
on the cardinality $\kappa$ of $H$. If $\kappa \leq \aleph_0$, then the
${\scriptstyle{{\bf LH}}}\mathfrak{Y}$-group $H$ is actually contained
in ${\scriptstyle{{\bf H}}}\mathfrak{Y}$ and we are done by the previous
discussion. If $\kappa$ is uncountable, then we may express $H$ as a
continuous ascending union of subgroups $(H_{\alpha})_{\alpha < \lambda}$
of cardinality $< \kappa$, which is indexed by a suitable ordinal $\lambda$.
The class $\mathfrak{Y}$ is subgroup-closed (cf.\ Proposition 4.1) and
hence this is also the case for the class ${\scriptstyle{\bf LH}}\mathfrak{Y}$;
in particular, $H_{\alpha} \in {\scriptstyle{\bf LH}}\mathfrak{Y}$ for all
$\alpha$. Invoking our induction hypothesis, we conclude that
$\mbox{res}_{H_{\alpha}}^GM$ is contained in ${\tt GInj}(kH_{\alpha})$
for all $\alpha$. Then, Lemma 4.3 below implies that the $kH$-module
$\mbox{res}_H^GM$ has Gorenstein injective dimension $\leq 1$. This is
the case for any $M \in \mathcal{I}(kG)$ and hence we can show as above
that $\mbox{res}_H^GM \in {\tt GInj}(kH)$.

(ii) We have to show that $\Phi_{inj}\mathfrak{Y} = \mathfrak{Y}$,
i.e.\ that $\Phi_{inj}\mathfrak{Y} \subseteq \mathfrak{Y}$. To that
end, we consider a $\Phi_{inj}\mathfrak{Y}$-group $G$ and prove that
$G \in \mathfrak{Y}$, by showing that
$\mathcal{I}(kG) \subseteq {\tt GInj}(kG)$. Let
$M \in \mathcal{I}(kG)$ and consider a short exact sequence of $kG$-modules
\begin{equation}
 0 \longrightarrow M \longrightarrow I
   \longrightarrow N \longrightarrow 0 ,
\end{equation}
where $I \in {\tt GInj}(kG)$ and $N \in \, \! ^{\perp}{\tt GInj}(kG)$. We
fix a $\mathfrak{Y}$-subgroup $H \subseteq G$. Since $M,I \in \mathcal{I}(kG)$,
it follows that $\mbox{res}_H^GM, \mbox{res}_H^GI \in \mathcal{I}(kH)$
and hence $\mbox{res}_H^GM, \mbox{res}_H^GI \in {\tt GInj}(kH)$. The class
of Gorenstein injective $kH$-modules is closed under cokernels of monomorphisms
and hence the $kH$-module $\mbox{res}_H^GN$ is also contained in ${\tt GInj}(kH)$.
On the other hand, coinduction from $H$ to $G$ maps ${\tt GInj}(kH)$ into
${\tt GInj}(kG)$ and hence the restriction-coinduction isomorphism
\[ \mbox{Ext}^1_{kH}(\mbox{res}_H^GN, \_\!\_ ) \simeq
   \mbox{Ext}^1_{kG}(N, \mbox{coind}_H^G \_\!\_) \]
shows that $\mbox{res}_H^GN \in \, \! ^{\perp}{\tt GInj}(kH)$. It follows
that $\mbox{res}_H^GN \in \, \! ^{\perp}{\tt GInj}(kH) \cap {\tt GInj}(kH)$
and hence the $kH$-module $\mbox{res}_H^GN$ is injective. Since this is
the case for any $\mathfrak{Y}$-subgroup $H \subseteq G$ and $G$ is a
$\Phi_{inj} \mathfrak{Y}$-group, we conclude that $\mbox{id}_{kG}N < \infty$;
in particular, the $kG$-module $N$ has finite Gorenstein injective dimension.
Then, the short exact sequence (5) shows that $M$ has finite Gorenstien
injective dimension as well. Having proved that any $kG$-module
$M \in \mathcal{I}(kG)$ has finite Gorenstein injective dimension, the closure
of $\mathcal{I}(kG)$ under direct products implies that there is an upper bound
on the Gorenstein injective dimensions of $\mathcal{I}(kG)$-modules. Arguing as
in (i), we can now show that any $\mathcal{I}(kG)$-module is Gorenstein injective.
\hfill $\Box$

\vspace{0.1in}

\noindent
We now state and prove the following variant of a result by Benson (cf.\
\cite[Lemma 4.5]{Ben}), that we used in the proof of Theorem 4.2(i).

\begin{Lemma}
Let $H$ be a group, which is expressed as a continuous ascending
union of subgroups $(H_{\alpha})_{\alpha < \lambda}$ indexed by
a suitable ordinal $\lambda$. If $N$ is a $kH$-module, which is
Gorenstein injective as a $kH_{\alpha}$-module for any
$\alpha < \lambda$, then $\mbox{Gid}_{kH}N \leq 1$.
\end{Lemma}
\vspace{-0.05in}
\noindent
{\em Proof.}
We fix a $kH$-module $M \in \, \! ^{\perp}{\tt GInj}(kH)$ and
consider the $kH$-module
$M_{\alpha} = \mbox{ind}_{H_{\alpha}}^H\mbox{res}_{H_{\alpha}}^HM$
for each $\alpha < \lambda$. Since coinduction from $H_{\alpha}$ to
$H$ maps ${\tt GInj}(kH_{\alpha})$ into ${\tt GInj}(kH)$, we may
conclude from the restriction-coinduction isomorphism (as in the
proof of Theorem 4.2(ii) above) that
$\mbox{res}_{H_{\alpha}}^HM \in \, \! ^{\perp}{\tt GInj}(kH_{\alpha})$
for all $\alpha < \lambda$. Hence, our assumption on $N$ implies that
\[ \mbox{Ext}_{kH}^i(M_{\alpha},N) = \mbox{Ext}_{kH_{\alpha}}^i
   \! \left( \mbox{res}_{H_{\alpha}}^HM , \mbox{res}_{H_{\alpha}}^HN
   \right) \! = 0 \]
for all $\alpha < \lambda$ and all $i \geq 1$. The inclusions
$H_{\alpha} \hookrightarrow H_{\beta}$, for $\alpha \leq \beta < \lambda$,
induce a continuous direct system of $kH$-modules
$(M_{\alpha})_{\alpha < \lambda}$ with surjective structure maps.
The short exact sequence
\[ 0 \longrightarrow K_{\alpha} \longrightarrow M_0
     \longrightarrow M_{\alpha} \longrightarrow 0 , \]
where $K_{\alpha}$ is the kernel of the structure map
$M_0 \longrightarrow M_{\alpha}$, is the $\alpha$-th term
of a continuous direct system of short exact sequences. The
colimit of the latter direct system of short exact sequences
is the short exact sequence of $kH$-modules
\begin{equation}
 0 \longrightarrow K \longrightarrow M_0
   \longrightarrow M \longrightarrow 0 .
\end{equation}
Here, $K$ is equal to the continuous ascending union of its
submodules $(K_{\alpha})_{\alpha < \lambda}$. Since
$K_{\alpha +1}/K_{\alpha}$ can be identified with the kernel
of the (surjective) structure map
$M_{\alpha} \longrightarrow M_{\alpha +1}$, we conclude that
$\mbox{Ext}^i_{kH}(K_{\alpha +1}/K_{\alpha},N) = 0$ for all
$\alpha < \lambda$ and all $i \geq 1$. Hence, Eklof's lemma
\cite[Lemma 6.2]{GT} implies that $\mbox{Ext}^i_{kH}(K,N)=0$
for all $i \geq 1$. Then, the short exact sequence (6) implies
that $\mbox{Ext}^i_{kH}(M,N)=0$ for all $i \geq 2$. Since the
latter equality is valid for any $kH$-module
$M \in \, \! ^{\perp}{\tt GInj}(kH)$, we conclude that
$\mbox{Gid}_{kH}N \leq 1$, as needed. \hfill $\Box$

\vspace{0.1in}

\noindent
The class $\mathfrak{Y}$ depends upon the ring $k$. In order for the
trivial group to be contained in $\mathfrak{Y}$, it is necessary that
$\mathcal{I}(k)={\tt GInj}(k)$. We note that the latter condition is
satisfied if, for example, $k$ is weakly Gorenstein regular, i.e.\ if
$\mbox{sfli} \, k < \infty$; this follows from \cite[Corollary 5.9]{Sto}.

\begin{Proposition}
If $\mathcal{I}(k) = {\tt GInj}(k)$, then all finite groups are
contained in $\mathfrak{Y}$.
\end{Proposition}
\vspace{-0.05in}
\noindent
{\em Proof.}
Assume that $\mathcal{I}(k) = {\tt GInj}(k)$ and let $G$ be a finite
group. We shall prove that $G \in \mathfrak{Y}$, i.e.\ that
$\mathcal{I}(kG) \subseteq {\tt GInj}(kG)$. To that end, we have to
show that any acyclic complex of injective $kG$-modules remains acyclic
upon applying the functor $\mbox{Hom}_{kG}(I,\_\!\_)$ for any injective
$kG$-module $I$. Since any injective $kG$-module is a direct summand of
a $kG$-module of the form $\mbox{coind}_1^GJ \simeq \mbox{ind}_1^GJ$,
where $J$ is an injective $k$-module, the result follows from the
induction-restriction isomorphism
\[ \mbox{Hom}_{kG} \! \left( \mbox{ind}_1^GJ,\_\!\_ \right) \! \simeq
   \mbox{Hom}_k \! \left( J,\mbox{res}_1^G\_\!\_ \right) \]
and our assumption that $\mathcal{I}(k) = {\tt GInj}(k)$. \hfill $\Box$

\vspace{0.1in}

\noindent
Let $\mathfrak{Y}_{fin}$ be the smallest class of groups, which
contains all finite groups and is ${\scriptstyle{{\bf LH}}}$-closed
and $\Phi_{inj}$-closed. Groups in class $\mathfrak{Y}_{fin}$ admit
a hierarchical description with base the class of finite groups, as
explained in $\S $1.IV.

\begin{Corollary}
If $\mathcal{I}(k) = {\tt GInj}(k)$, then all
$\mathfrak{Y}_{fin}$-groups are contained in $\mathfrak{Y}$.
\end{Corollary}
\vspace{-0.05in}
\noindent
{\em Proof.}
Proposition 4.4 and Theorem 4.2 imply that
$\mathfrak{Y}_{fin} \subseteq \mathfrak{Y}$. \hfill $\Box$

\vspace{0.1in}

\noindent
{\bf Remarks 4.6.}
(i) Let $k$ be a commutative ring. As we remarked above,
the equality $\mathcal{I}(k) = {\tt GInj}(k)$ is satisfied
if $k$ is weakly Gorenstein regular. It follows that the pairs
$(k,G)$ considered in Corollary 4.5 include those considered
in \cite{Bi2}, where it was assumed that $k$ has finite global
dimension and $G$ is either an
${\scriptstyle{{\bf LH}}}\mathfrak{F}$-group or a group of
type $\Phi$. (Indeed, if $\mbox{gl.dim} \, k < \infty$ and
$G$ is a group of type $\Phi$, then $kG$ is Gorenstein
regular. Invoking Proposition 1.6, we conclude that $G$ is
contained in $\Phi_{inj}\mathfrak{F}$ and hence
in $\mathfrak{Y}_{fin}$.)

(ii) Assume that $k = \mathbb{Z}$ is the ring of integers.
It would be of interest to know of a group $G$ which is
contained in $\mathfrak{Y}_{fin}$, but is not an
${\scriptstyle{{\bf LH}}}\mathfrak{F}$-group nor a group of
type $\Phi$.
\addtocounter{Lemma}{1}

\section{Gorenstein flat modules over group algebras and duality}

\noindent
In this section, we consider groups $G$, over which the $kG$-modules
whose Pontryagin dual is Gorenstein injective are necessarily
Gorenstein flat. Since the dual $DM$ of any Gorenstein flat
$kG$-module $M$ is Gorenstein injective, we have an inclusion
$D{\tt GFlat}(kG) \subseteq {\tt GInj}(kG)$, i.e.\
\[ {\tt GFlat}(kG) \subseteq D^{-1}{\tt GInj}(kG) . \]
Having fixed $k$, we consider the group class
$\mathfrak{Z} = \mathfrak{Z}(k)$, which is defined by letting
\[ \mathfrak{Z} = \mathfrak{Y} \cap
   \{ G : D^{-1}{\tt GInj}(kG) \subseteq {\tt GFlat}(kG) \} =
   \mathfrak{Y} \cap
   \{ G : {\tt GFlat}(kG) = D^{-1}{\tt GInj}(kG) \} . \]
For groups in class $\mathfrak{Z}$, all of the general questions
considered in the Introduction, regarding the behaviour of the
associated group algebra with respect to Gorenstein homological
algebra notions, admit a positive answer.

\begin{Proposition}
The class $\mathfrak{Z}$ is subgroup-closed and for any
$\mathfrak{Z}$-group $G$ the following hold:

(i) The cotorsion pair
$\left( {\tt GProj}(kG) , {\tt GProj}(kG)^{\perp} \right)$ is
complete.

(ii) ${\tt GProj}(kG) \subseteq {\tt GFlat}(kG)$.

(iii) ${\tt GFlat}(kG) = D^{-1}{\tt GInj}(kG)$.

(iv) Any acyclic complex of projective $kG$-modules is totally
acyclic, i.e.\ it remains acyclic after applying the functor
$\mbox{Hom}_{kG}(\_\!\_,P)$ for any projective $kG$-module $P$.

(v) Any acyclic complex of flat $kG$-modules is totally acyclic,
i.e.\ it remains acyclic after applying the functor
$I \otimes_{kG} \_\!\_$ for any injective $kG$-module $I$.

(vi) Any acyclic complex of injective $kG$-modules is totally
acyclic, i.e.\ it remains acyclic after applying the functor
$\mbox{Hom}_{kG}(I,\_\!\_)$ for any injective $kG$-module $I$.
\end{Proposition}
\vspace{-0.05in}
\noindent
{\em Proof.}
Assertion (iii) is part of the definition of $\mathfrak{Z}$.
Since $\mathfrak{Z} \subseteq \mathfrak{Y}$, assertions (i),
(ii), (iv), (v) and (vi) follow from Proposition 4.1.

It remains to show that the class $\mathfrak{Z}$ is subgroup-closed.
To that end, we consider a $\mathfrak{Z}$-group $G$ and let
$H \subseteq G$ be a subgroup. Then, $G$ is a $\mathfrak{Y}$-group
and $D^{-1}{\tt GInj}(kG) \subseteq {\tt GFlat}(kG)$. We have to
show that $H$ is a $\mathfrak{Z}$-group. Since $\mathfrak{Y}$ is
subgroup-closed (cf.\ Proposition 4.1), $H$ is a $\mathfrak{Y}$-group
as well. In order to show that
$D^{-1}{\tt GInj}(kH) \subseteq {\tt GFlat}(kH)$, we consider
a $kH$-module $M$, for which $DM \in {\tt GInj}(kH)$. Then,
$D \, \mbox{ind}_H^GM = \mbox{coind}_H^GDM \in {\tt GInj}(kG)$
and hence our assumption on $G$ implies that
$\mbox{ind}_H^GM \in {\tt GFlat}(kG)$. Since
$H \in \mathfrak{Y} \subseteq \mathfrak{X}$ and restriction maps
$\mathcal{F}(kG)$ into $\mathcal{F}(kH) = {\tt GFlat}(kH)$, it
follows that $\mbox{res}_H^G\mbox{ind}_H^GM \in {\tt GFlat}(kH)$.
The $kH$-module $M$ being a direct summand of
$\mbox{res}_H^G\mbox{ind}_H^GM$, we conclude that
$M \in {\tt GFlat}(kH)$, as needed. \hfill $\Box$

\vspace{0.1in}

\noindent
In order to provide examples of $\mathfrak{Z}$-groups, we shall
prove below that $\mathfrak{Z}$ is also closed under the operations
${\scriptstyle{{\bf LH}}}$ and $\Phi_{inj}$. To that end, we shall
use Bouchiba's criterion (cf.\ Proposition 1.1): A $kG$-module $M$
is Gorenstein flat, if (and only if) $M$ has finite Gorenstein
flat dimension and the dual $DM$ of $M$ is Gorenstein injective.

\begin{Theorem}
The class $\mathfrak{Z}$ is (i) ${\scriptstyle{{\bf LH}}}$-closed
and (ii) $\Phi_{inj}$-closed.
\end{Theorem}
\vspace{-0.05in}
\noindent
{\em Proof.}
(i) We have to show that
${\scriptstyle{{\bf LH}}}\mathfrak{Z} = \mathfrak{Z}$, i.e.\
that
${\scriptstyle{{\bf LH}}}\mathfrak{Z} \subseteq \mathfrak{Z}$.
To that end, let us consider a group
$G \in {\scriptstyle{{\bf LH}}}\mathfrak{Z}$. Since
$\mathfrak{Z} \subseteq \mathfrak{Y}$, we have
${\scriptstyle{{\bf LH}}}\mathfrak{Z} \subseteq
 {\scriptstyle{{\bf LH}}}\mathfrak{Y} = \mathfrak{Y}$
and hence $G \in \mathfrak{Y}$. It follows that
$\mathcal{I}(kG) = {\tt GInj}(kG)$. Since $\mathfrak{Y}$
is subgroup-closed, any subgroup $H \subseteq G$ is also
a $\mathfrak{Y}$-group and hence
$\mathcal{I}(kH) = {\tt GInj}(kH)$. In order to show that
$G$ is a $\mathfrak{Z}$-group, it remains to show that
$D^{-1}{\tt GInj}(kG) \subseteq {\tt GFlat}(kG)$.

First of all, we shall use the standard argument, in order
to show that $\mbox{res}_H^GM \in {\tt GFlat}(kH)$ for all
${\scriptstyle{\bf H}}\mathfrak{Z}$-subgroups
$H \subseteq G$ and all $kG$-modules $M$ with
$DM \in \mathcal{I}(kG) = {\tt GInj}(kG)$. We use transfinite
induction on the ordinal $\alpha$, which is such that
$H \in {\scriptstyle{{\bf H}}}_{\alpha}\mathfrak{Z}$. If
$\alpha = 0$, then $H \in \mathfrak{Z}$. Since
$D\mbox{res}_H^GM = \mbox{res}_H^GDM \in \mathcal{I}(kH) =
 {\tt GInj}(kH)$,
we conclude that $\mbox{res}_H^GM \in {\tt GFlat}(kH)$. We
now assume that $\alpha > 0$ and the result is true for all
$kG$-modules $M$ with $DM \in \mathcal{I}(kG) = {\tt GInj}(kG)$
and all ${\scriptstyle{{\bf H}}}_{\beta}\mathfrak{Z}$-subgroups
$H \subseteq G$ for all $\beta < \alpha$. Let us consider a
$kG$-module $M$ with $DM \in \mathcal{I}(kG) = {\tt GInj}(kG)$
and let $H$ be an
${\scriptstyle{{\bf H}}}_{\alpha}\mathfrak{Z}$-subgroup of
$G$. Then, there exists an exact sequence of $kH$-modules
\[ 0 \longrightarrow M_r \longrightarrow \cdots
     \longrightarrow M_1 \longrightarrow M_0
     \longrightarrow \mbox{res}_H^GM \longrightarrow 0 , \]
where each $M_i$ is a direct sum of $kH$-modules of the form
$\mbox{ind}_F^H \mbox{res}_F^GM$, for a suitable
${\scriptstyle{{\bf H}}}_{\beta}\mathfrak{Z}$-subgroup
$F \subseteq H$ for some $\beta < \alpha$.
Since induction from $F$ to $H$ maps the class ${\tt GFlat}(kF)$
into ${\tt GFlat}(kH)$ for any subgroup $F \subseteq H$, we
conclude that $M_i \in {\tt GFlat}(kH)$ for all $i=0,1, \ldots ,r$.
It follows that the $kH$-module $\mbox{res}_H^GM$ has Gorenstein
flat dimension $\leq r$. Since
$D\mbox{res}_H^GM = \mbox{res}_H^GDM \in \mathcal{I}(kH) =
 {\tt GInj}(kH)$,
Bouchiba's criterion implies that
$\mbox{res}_H^GM \in {\tt GFlat}(kH)$.

We shall now prove that $\mbox{res}_H^GM \in {\tt GFlat}(kH)$ for
all ${\scriptstyle{\bf LH}}\mathfrak{Z}$-subgroups $H \subseteq G$
and all $kG$-modules $M$ with $DM \in \mathcal{I}(kG)={\tt GInj}(kG)$.
(Then, letting $H=G$, we will conclude that
$D^{-1}{\tt GInj}(kG) \subseteq {\tt GFlat}(kG)$.) We proceed by
induction on the cardinality $\kappa$ of $H$. If $\kappa \leq \aleph_0$,
then the ${\scriptstyle{{\bf LH}}}\mathfrak{Z}$-group $H$ is actually
contained in ${\scriptstyle{{\bf H}}}\mathfrak{Z}$ and we are done by
the previous discussion. If $\kappa$ is uncountable, we may express
$H$ as a continuous ascending union of subgroups
$(H_{\alpha})_{\alpha < \lambda}$ of cardinality $< \kappa$, which
is indexed by a suitable ordinal $\lambda$. The class $\mathfrak{Z}$
is subgroup-closed (cf.\ Proposition 5.1) and hence this is also the
case for the class ${\scriptstyle{\bf LH}}\mathfrak{Z}$; in particular,
$H_{\alpha} \in {\scriptstyle{\bf LH}}\mathfrak{Z}$ for all $\alpha$.
Let $M$ be a $kG$-module for which $DM \in \mathcal{I}(kG)={\tt GInj}(kG)$.
Our induction hypothesis implies that
$\mbox{res}_{H_{\alpha}}^GM \in {\tt GFlat}(kH_{\alpha})$ and hence
$\mbox{ind}_{H_{\alpha}}^H\mbox{res}_{H_{\alpha}}^GM \in {\tt GFlat}(kH)$
for all $\alpha < \lambda$. As $\mbox{res}_H^GM$ is the filtered
colimit of the system
$\left( \mbox{ind}_{H_{\alpha}}^H\mbox{res}_{H_{\alpha}}^GM \right) \!
 _{\alpha < \lambda}$,
an application of \cite[Corollary 4.12]{SS} implies that
$\mbox{res}_H^GM \in {\tt GFlat}(kH)$, as needed.

(ii) We have to show that $\Phi_{inj} \mathfrak{Z} = \mathfrak{Z}$,
i.e.\ that $\Phi_{inj} \mathfrak{Z} \subseteq \mathfrak{Z}$. To
that end, consider a group $G \in \Phi_{inj} \mathfrak{Z}$. Since
$\mathfrak{Z} \subseteq \mathfrak{Y}$, we also have
$\Phi_{inj}\mathfrak{Z} \subseteq
 \Phi_{inj}\mathfrak{Y} = \mathfrak{Y}$
and hence $G \in \mathfrak{Y}$; therefore,
$\mathcal{I}(kG) = {\tt GInj}(kG)$. In order to show that $G$
is a $\mathfrak{Z}$-group, it remains to prove the inclusion
$D^{-1}{\tt GInj}(kG) \subseteq {\tt GFlat}(kG)$.

Let $M$ be a $kG$-module, such that
$DM \in \mathcal{I}(kG) = {\tt GInj}(kG)$, and consider a short
exact sequence of $kG$-modules
\begin{equation}
 0 \longrightarrow N \longrightarrow F
   \longrightarrow M \longrightarrow 0 ,
\end{equation}
where $N \in {\tt GFlat}(kG)^{\perp}$ and
$F \in {\tt GFlat}(kG)$. We also fix a $\mathfrak{Z}$-subgroup
$H \subseteq G$. Since the $kG$-modules $DM$ and $DF$ are
contained in $\mathcal{I}(kG) = {\tt GInj}(kG)$, the
$kH$-modules $D\mbox{res}_H^GM = \mbox{res}_H^GDM$ and
$D\mbox{res}_H^GF = \mbox{res}_H^GDF$ are contained in
$\mathcal{I}(kH) = {\tt GInj}(kH)$. Hence, we conclude that
$\mbox{res}_H^GM, \mbox{res}_H^GF \in {\tt GFlat}(kH)$. The
closure of ${\tt GFlat}(kH)$ under kernels of epimorphisms
implies that $\mbox{res}_H^GN \in {\tt GFlat}(kH)$ as well.
On the other hand, induction from $H$ to $G$ maps
${\tt GFlat}(kH)$ into ${\tt GFlat}(kG)$ and hence the
induction-restriction isomorphism
\[ \mbox{Ext}^1_{kG}(\mbox{ind}_H^G \_\!\_ , N) \simeq
   \mbox{Ext}^1_{kH}(\_\!\_ , \mbox{res}_H^GN) \]
shows that $\mbox{res}_H^GN \in {\tt GFlat}(kH)^{\perp}$.
Then,
$\mbox{res}_H^GN \in {\tt GFlat}(kH) \cap {\tt GFlat}(kH)^{\perp}$
and hence the $kH$-module $\mbox{res}_H^GN$ is flat-cotorsion.
Since this is the case for any $\mathfrak{Z}$-subgroup
$H \subseteq G$ and
$G \in \Phi_{inj} \mathfrak{Z} \subseteq
 \Phi_{flat} \mathfrak{Z}$
(cf.\ Lemma 1.5), we conclude that $\mbox{fd}_{kG}N < \infty$;
in particular, the $kG$-module $N$ has finite Gorenstein flat
dimension. The short exact sequence (7) shows that $M$ has finite
Gorenstein flat dimension as well. Since $DM$ is Gorenstein injective,
Bouchiba's criterion implies that $M$ is Gorenstein flat. \hfill $\Box$

\vspace{0.1in}

\noindent
The class $\mathfrak{Z}$ depends upon the ring $k$. In order
for the trivial group to be contained therein, it is necessary
that $\mathcal{I}(k) = {\tt GInj}(k)$ and
$D^{-1}{\tt GInj}(k) = {\tt GFlat}(k)$. As we have already
noted, \cite[Corollary 5.9]{Sto} implies that the equality
$\mathcal{I}(k)={\tt GInj}(k)$ holds if $k$ is weakly Gorenstein
regular, i.e.\ if $\mbox{sfli} \, k < \infty$. On the other hand,
an application of \cite[Theorem 3.6]{H} shows that the equality
$D^{-1}{\tt GInj}(k) = {\tt GFlat}(k)$ holds if $k$ is coherent.

\begin{Proposition}
If $\mathcal{I}(k) = {\tt GInj}(k)$ and
$D^{-1}{\tt GInj}(k)={\tt GFlat}(k)$, then all finite groups are
contained in $\mathfrak{Z}$.
\end{Proposition}
\vspace{-0.05in}
\noindent
{\em Proof.}
Let $G$ be a finite group. Since $\mathcal{I}(k) = {\tt GInj}(k)$,
Proposition 4.4 implies that $G \in \mathfrak{Y}$. Assuming that
we also have $D^{-1}{\tt GInj}(k)={\tt GFlat}(k)$, we shall prove
that $D^{-1}{\tt GInj}(kG) \subseteq {\tt GFlat}(kG)$. To that end,
let $M$ be a $kG$-module, such that $DM$ is Gorenstein injective.
Since
$D\mbox{res}_1^GM = \mbox{res}_1^GDM \in \mathcal{I}(k) =
 {\tt GInj}(k)$,
the additional assumption made on $k$ implies that
$\mbox{res}_1^GM \in {\tt GFlat}(k)$. Then, the induced $kG$-module
$\mbox{ind}_1^G\mbox{res}_1^GM$ is Gorenstein flat and hence the
$kG$-module $M$ is itself Gorenstein flat, in view of Proposition
1.2(iii). \hfill $\Box$

\vspace{0.1in}

\noindent
We recall that $\mathfrak{Y}_{fin}$ was defined as the smallest
class of groups, which contains all finite groups and is
${\scriptstyle{{\bf LH}}}$-closed and $\Phi_{inj}$-closed. The
following result complements Corollary 4.5.

\begin{Corollary}
Assume that $\mathcal{I}(k) = {\tt GInj}(k)$,
$D^{-1}{\tt GInj}(k) = {\tt GFlat}(k)$ and
$G \in \mathfrak{Y}_{fin}$. Then:

(i) The cotorsion pair
$\left( {\tt GProj}(kG) , {\tt GProj}(kG)^{\perp} \right)$ is
complete.

(ii) ${\tt GProj}(kG) \subseteq {\tt GFlat}(kG)$.

(iii) ${\tt GFlat}(kG) = D^{-1}{\tt GInj}(kG)$.

(iv) Any acyclic complex of projective $kG$-modules is totally
acyclic, i.e.\ it remains acyclic after applying the functor
$\mbox{Hom}_{kG}(\_\!\_,P)$ for any projective $kG$-module $P$.

(v) Any acyclic complex of flat $kG$-modules is totally acyclic,
i.e.\ it remains acyclic after applying the functor
$I \otimes_{kG} \_\!\_$ for any injective $kG$-module $I$.

(vi) Any acyclic complex of injective $kG$-modules is totally
acyclic, i.e.\ it remains acyclic after applying the functor
$\mbox{Hom}_{kG}(I,\_\!\_)$ for any injective $kG$-module $I$.
\end{Corollary}
\vspace{-0.05in}
\noindent
{\em Proof.}
Proposition 5.3 shows that $\mathfrak{Z}$ contains the class
of finite groups. Then, Theorem 5.2 implies that we also have
$\mathfrak{Y}_{fin} \subseteq \mathfrak{Z}$. Hence, the result
follows from Proposition 5.1.\hfill $\Box$

\vspace{0.1in}

\noindent
We recall that the group algebra $kG$ is Gorenstein regular if
$kG$ has finite Gorenstein global dimension; this is the case
if and only if both invariants $\mbox{silp} \, kG$ and
$\mbox{spli} \, kG$ are finite. A weaker condition than Gorenstein
regularity is to require that $kG$ is weakly Gorenstein regular,
i.e.\ that the Gorenstein weak global dimension of $kG$ be finite;
the latter condition is equivalent to the finiteness of the invariant
$\mbox{sfli} \, kG$.

\begin{Proposition}
If $kG$ is weakly Gorenstein regular, then $G \in \mathfrak{Z}$.
\end{Proposition}
\vspace{-0.05in}
\noindent
{\em Proof.}
Since $\mbox{sfli} \, kG$ is finite, the functor
$\mbox{Hom}_{kG}(I,\_\!\_)$ preserves the acyclicity of any acyclic
complex of injective $kG$-modules for any injective $kG$-module $I$;
cf.\ \cite[Corollary 5.9]{Sto}. Hence, it follows that
$\mathcal{I}(kG)={\tt GInj}(kG)$, so that $G \in \mathfrak{Y}$.
It remains to show that $D^{-1}{\tt GInj}(kG) = {\tt GFlat}(kG)$.
In other words, we have to show that any $kG$-module $M$ whose dual
$DM$ is Gorenstein injective is necessarily Gorenstein flat. Since
$M$ has finite Gorenstein flat dimension, the latter claim follows
from Bouchiba's criterion. \hfill $\Box$

\vspace{0.1in}

\noindent
Having fixed the commutative ring $k$, let $\mathfrak{Z}_{Gor}$
be the smallest class of groups which contains all groups $G$ for
which $kG$ is weakly Gorenstein regular\footnote{In particular,
$\mathfrak{Z}_{Gor}$ contains all groups $G$ for which $kG$ is
Gorenstein regular.} and is ${\scriptstyle{{\bf LH}}}$-closed
and $\Phi_{inj}$-closed. Groups in class $\mathfrak{Z}_{Gor}$
admit a hierarchical description, as explained in $\S $1.IV.

\begin{Corollary}
Let $G$ be a $\mathfrak{Z}_{Gor}$-group. Then:

(i) The cotorsion pair
$\left( {\tt GProj}(kG) , {\tt GProj}(kG)^{\perp} \right)$ is
complete.

(ii) ${\tt GProj}(kG) \subseteq {\tt GFlat}(kG)$.

(iii) ${\tt GFlat}(kG) = D^{-1}{\tt GInj}(kG)$.

(iv) Any acyclic complex of projective $kG$-modules is totally
acyclic, i.e.\ it remains acyclic after applying the functor
$\mbox{Hom}_{kG}(\_\!\_,P)$ for any projective $kG$-module $P$.

(v) Any acyclic complex of flat $kG$-modules is totally acyclic,
i.e.\ it remains acyclic after applying the functor
$I \otimes_{kG} \_\!\_$ for any injective $kG$-module $I$.

(vi) Any acyclic complex of injective $kG$-modules is totally
acyclic, i.e.\ it remains acyclic after applying the functor
$\mbox{Hom}_{kG}(I,\_\!\_)$ for any injective $kG$-module $I$.
\end{Corollary}
\vspace{-0.05in}
\noindent
{\em Proof.}
Proposition 5.5 and Theorem 5.2 imply that
$\mathfrak{Z}_{Gor} \subseteq \mathfrak{Z}$. Hence, the result
follows from Proposition 5.1. \hfill $\Box$

\vspace{0.1in}

\noindent
As a variation of the definition of $\mathfrak{Z}_{Gor}$, let
$\mathfrak{X}_{Gor}$ be the smallest class of groups which
contains all groups $G$ for which $kG$ is weakly Gorenstein
regular and is ${\scriptstyle{{\bf LH}}}$-closed,
$\Phi_{proj}$-closed and $\Phi_{flat}$-closed. Groups in class
$\mathfrak{X}_{Gor}$ admit a hierarchical description, as
explained in $\S $1.IV.

\begin{Corollary}
Any $\mathfrak{X}_{Gor}$-group is contained in $\mathfrak{X}$.
\end{Corollary}
\vspace{-0.05in}
\noindent
{\em Proof.}
Since $\mathfrak{Z} \subseteq \mathfrak{Y} \subseteq \mathfrak{X}$,
Proposition 5.5 implies that any group $G$ for which $kG$ is weakly
Gorenstein regular is necessarily an $\mathfrak{X}$-group. The
inclusion $\mathfrak{X}_{Gor} \subseteq \mathfrak{X}$ is therefore
an immediate consequence of Theorem 3.3. \hfill $\Box$

\vspace{0.1in}

\noindent
{\bf Remarks 5.8.}
(i) The pairs $(k,G)$ considered in Corollaries 5.6 and 5.7
(in fact, the pairs considered in Proposition  5.5) include
those considered in \cite{R, RY}, where it was assumed that
$k$ is (weakly) Gorenstein regular and $G$ has finite Gorenstein
(co-)homological dimension.

(ii) Assume that $k = \mathbb{Z}$ is the ring of integers. It
would be of interest to know of a group $G$ which is contained
in $\mathfrak{Y}_{Gor}$ or $\mathfrak{X}_{Gor}$, but is not an
${\scriptstyle{{\bf LH}}}\mathfrak{F}$-group nor a group of type
$\Phi$.
\addtocounter{Lemma}{1}

\vspace{0.1in}

\noindent
{\small {\sc Department of Mathematics,
             University of Athens,
             Athens 15784,
             Greece}}

\noindent
{\em E-mail addresses:} {\tt emmanoui@math.uoa.gr} and
                        {\tt otalelli@math.uoa.gr}

\end{document}